\documentclass[centertags,reqno,a4paper]{amsart}
\usepackage{latexsym,amsmath,amssymb,amscd,amsthm,epsfig}
\usepackage[all]{xy}
\newcommand{\N}{{\mathbb N}}
\newcommand{\A}{{\mathbb A}}
\newcommand{\Z}{{\mathbb Z}}
\newcommand{\R}{{\mathbb R}}
\newcommand{\C}{{\mathbb C}}

\renewcommand{\P}{{\mathbb P}}

\newcommand{\G}{\mathbb{G}}
\newcommand{\K}{\C}   
\newcommand{\M}{\mathcal{M}}

\DeclareMathOperator{\mt}{mt}
\DeclareMathOperator{\clg}{clg}
\DeclareMathOperator{\Cl}{\kc\!\!\;\ell\!\:}

\DeclareMathOperator{\tHilb}{Hilb\!\:}
\DeclareMathOperator{\tHil}{Hilb\!\:}

\newcommand{\lra}{\longrightarrow}

\newcommand{\Hilb}{\kh ilb\;\!}

\newcommand{\kc}{{\mathcal C}}

\newcommand{\kh}{{\mathcal H}}
\newcommand{\ki}{{\mathcal I}}
\newcommand{\kj}{{\mathcal J}}
\newcommand{\kk}{{K}}

\newcommand{\ko}{{\mathcal O}}

\newcommand{\ku}{{\mathcal U}}

\newcommand{\kx}{{\mathcal H}}

\newcommand{\kxf}{{\mathcal H}_0}
\newcommand{\fm}{\mathfrak{m}}
\newcommand{\fx}{\mathfrak{X}}

\newcommand{\fu}{\mathfrak{U}}

\newtheorem{lemma}{Lemma}[section]
\newtheorem{lemma and definition}[lemma]{Lemma and Definition}
\newtheorem{proposition}[lemma]{Proposition}
\newtheorem{theorem}{Theorem}
\newtheorem{corollary}[lemma]{Corollary}
\newtheorem{conjecture}[lemma]{Conjecture}

\theoremstyle{remark}
\newtheorem{remark}[lemma]{Remark}
\newtheorem{examples}{Examples$\!\!$} 
 
\newtheorem{definition}{Definition$\!$} 
\newtheorem{remdef}[lemma]{Remark and Definition} 
\newcommand{\Ies}{I^{\text{es}}}
\newcommand{\Iea}{I^{\text{ea}}}
\newcommand{\Iesf}{I^{\text{es}}_{\text{fix}}}
\newcommand{\Ieaf}{I^{\text{ea}}_{\text{fix}}}
\newcommand{\Xes}{X^{\text{es}}}
\newcommand{\Xea}{X^{\text{ea}}}
\newcommand{\Xesf}{X^{\text{es}}_{\text{fix}}}
\newcommand{\Xeaf}{X^{\text{ea}}_{\text{fix}}}
\newcommand{\Vd}{V_d^{\text{\it irr}}(S_1,\dots,S_r)}
\newcommand{\Vdi}{V_d^{\text{\it irr}}}

\newcommand{\Vr}{V_{\text{\it reg}}}
\newcommand{\Vg}{V_{\text{\it gen}}}

\begin{document}

\title{Castelnuovo function,
  zero-dimensional  schemes and singular plane curves}
\author{
Gert-Martin Greuel}
\address{Universit\"at Kaiserslautern\\
Fachbereich Mathematik\\
Erwin-Schr\"odinger-Stra\ss e\\
D -- 67663 Kaiserslautern\\
e-mail: greuel@mathematik.uni-kl.de}
\author{
Christoph Lossen}
\address{
Universit\"at Kaiserslautern\\
Fachbereich Mathematik\\
Erwin-Schr\"odinger-Stra\ss e\\
D -- 67663 Kaiserslautern\\
e-mail: lossen@mathematik.uni-kl.de}
\author{ 
Eugenii Shustin}
\address{Tel Aviv University\\
School of Mathematical Sciences\\
Ramat Aviv\\
ISR -- Tel Aviv 69978\\
e-mail: shustin@math.tau.ac.il
}
\thanks{\noindent Work on this paper has been partially supported by the 
Hermann Minkowski -- Minerva Center for Geometry at Tel Aviv University
and Grant No.\ G 039-304.01/95 from the
German Israeli Foundation for Research and Development.}

\begin{abstract}
We study families $V$ of curves in $\P^2(\C)$ of degree $d$ having exactly $r$ 
singular points of given topological or analytic types.
We derive new sufficient conditions for $V$ to be T-smooth (smooth of 
the expected dimension), respectively to be irreducible. 
For T-smoothness these conditions involve new invariants 
of curve singularities and are conjectured to be asymptotically proper, i.e., 
optimal up to a constant factor. To obtain the results, we study the 
Castelnuovo function, prove the irreducibility of the Hilbert scheme of 
zero-dimensional schemes associated to a cluster of infinitely near points 
of the singularities and deduce new vanishing theorems for ideal sheaves of
zero-dimensional schemes in $\P^2$.
Moreover, we give a series of examples of cuspidal curves where the family 
$V$ is reducible, but where \mbox{$\pi_1(\P^2\!\setminus\!\!\:C)$} 
coincides (and is abelian) for all \mbox{$C\in V$}. 
\end{abstract}
\setcounter{tocdepth}{1}
\maketitle
\tableofcontents
\thispagestyle{empty}
\date{}

\section*{Introduction}
\subsection*{Statement of the problem and asymptotically proper bounds}
Singular algebraic curves, their existence, deformation, families (from the
local and global point of view) attract continuous attention of algebraic
geometers since the last century. The geometry of equisingular families 
of algebraic curves on smooth algebraic surfaces has been founded in
basic works of Pl\"ucker, Severi, Segre, Zariski, and has tight links and
finds important applications in singularity theory, topology of complex
algebraic curves and surfaces, and in real algebraic geometry.

In the present paper we consider the family \mbox{$\Vd$} of reduced 
irreducible complex 
plane curves of degree $d$ with $r$ isolated singular points of
given topological, or analytic types \mbox{$S_1,\dots,S_r$}
(further referred to as equisingular families, or ESF). 
The questions about the non-emptiness, smoothness, irreducibility and
dimension are basic in the geometry of ESF.
Except for the case of nodal
curves, no complete answers are known and one can hardly expect them.

Our goal, however, is to obtain {\it asymptotically proper}
sufficient conditions for ESF to have ``good" properties like being non-empty,
or smooth, or irreducible. The conditions should be expressed
in the form of bounds to numerical invariants of curves and
singularities such that, for the ``good'' properties to hold, the necessary
respectively sufficient conditions should be given by inequalities with the
same invariants but, maybe, with different absolute constants.
As an example, we mention our sufficient condition for the non-emptiness of
$\Vd$ with topological singularities $S_1,\dots,S_r$ \cite{GLS96,Lo,Lo1}
\begin{equation}
\label{****}
\textstyle{\sum\limits_{i=1}^r} \mu(S_i)<\tfrac{1}{46}\cdot (d+2)^2\ ,
\end{equation}
whereas the classically known necessary condition is
\begin{equation}
\label{??}
\textstyle{\sum\limits_{i=1}^r} \mu(S_i)\le(d-1)^2\ .
\end{equation}
In the present paper we obtain two qualitatively new bounds: one for the
smoothness of $\Vd$ and one
for the irreducibility. In particular, we show that the inequality
\begin{equation}
\textstyle{\sum\limits_{i=1}^r} \gamma(S_i)< d^2+6d+8\,,  \label{i1}
\end{equation}
where $\gamma(S)$ is a new singularity invariant (defined in Section
\ref{sec:smoothness es}), 
is sufficient for the smoothness and expected dimension 
(also called {\em T-property\/}) of $\Vd$.  
We expect (\ref{i1}) to be  asymptotically proper for  topological
singularities in the following sense: 

\begin{conjecture}\label{con1}   
There exists an absolute constant \mbox{$A>0$} such that for any 
topological singularity $S$ there are infinitely many pairs
\mbox{$(r,d)\in\N\!\:^2$}  such that
\mbox{$\Vdi(r\!\!\:\cdot\!\!\: S)$} is empty or not smooth or has dimension 
greater than the expected one and 
$$r\cdot\gamma(S)\:\le\: A\cdot d^2.$$
\end{conjecture}

\noindent
We know that the exponent $2$ of $d$ in the right-hand side of (\ref{i1})
cannot be raised in any reasonable sufficient criterion for T-property with the
left-hand side being the sum of local singularity invariants. Hence, for an
asymptotically proper sufficient criterion for T-property  the right-hand side
is correct. On the other hand, for the left-hand side of such a sufficient
criterion different invariants can be used. What we conjecture is that the new
invariant $\gamma(S)$ is the ``correct'' one for an asymptotically proper 
bound in the case of topological singularities.

The conjecture is known to be true for an infinite series of singularities of
types $A$ and $D$ (cf.\ \cite{Shu97,GLS97}) and it holds for  
ordinary singularities, because here the inequality (\ref{i1}) is implied by
\begin{equation}
4\cdot \#(\text{nodes\/})+18
\cdot \#(\text{triple points\/})+
\textstyle{\sum\limits_{\mt S_i>3}} 
\tfrac{16}{7} \cdot(\mt S_i)^2 \:<\:d^2+6d+8\,,\label{i3}
\end{equation}
(cf.\ Corollary \ref{Smoothn ord cor}) whereas the inequality
$$\textstyle{\sum\limits_{i=1}^r} \mt S_i(\mt S_i-1) \: \leq \: (d-1)(d-2)$$
is necessary for the existence of an irreducible curve with ordinary
singularities \mbox{$S_1,\dots,S_r$}.

\subsection*{New criteria for smoothness and irreducibility of equisingular
  families} 
We show that under condition (\ref{i1})
(with singularity invariants \mbox{$\gamma(S)\leq 
(\tau'(S)\!\!\;+\!\!\;1)^2$},
where $\tau'$ stands for the Tjurina number $\tau$ if $S$ is an analytic type
and for \mbox{$\tau^{\text{es}}=\mu-\text{modality}$} if $S$ is a topological
type) the family \mbox{$V=\Vd$} is either empty or smooth of
the expected dimension (Theorem \ref{Smoothness Theorem} in Section
\ref{sec:smoothness}).
In addition, for any curve
\mbox{$C\in \Vd$} the inequality
(\ref{i1}) with analytic invariants $\gamma$ is sufficient for the
independence of versal deformations of all singular points 
when varying in the space of plane curves of degree $d$.

This improves the previously known condition (cf.\ \cite{GLS97})
$$\textstyle{\sum\limits_{i=1}^r} (\tau'(S_i)+1)^2\:<\:d^2\, ,$$
mainly with respect to the singularity invariants in the left-hand side. For
instance, for an ordinary singular point $S$ of multiplicity $m$, considered up
to topological equivalence,
$$(\tau'(S)+1)^2=\bigl(\tfrac{m(m+1)}{2}-1\bigr)^2\ \sim\
\tfrac{1}{4} m^4,$$ 
whereas the invariant in the left-hand side of the new condition is
\mbox{$\gamma(S) \leq \frac{16}{7}\!\;m^2$} (cf.\ (\ref{i3})).

Another new result concerns the irreducibility of ESF. It says 
that under the conditions \mbox{$\max_i
\tau'(S_i)\leq \frac{2}{5}\!\;d-1$} and
\begin{equation}
\tfrac{25}{2}\cdot \#(\text{nodes})+18\cdot \#(
\text{cusps})+
\tfrac{10}{9} \cdot \!\textstyle{\sum\limits_{\tau'(S_i)\ge 3}}
 (\tau'(S_i)\!+\!2)^2 \,<\: d^2
\label{i2}
\end{equation}
the family $\Vd$ is irreducible (cf.\ Theorem \ref{Theorem 3} in Section
\ref{sec:irred}  with a slightly stronger statement).
The irreducibility criterion (\ref{i2}) improves the bounds
known before
\begin{equation}
\textstyle{\sum\limits_{i=1}^r} 
\mu(S_i) \,<\, \min\limits_{1\le i\le r} f(S_i) \cdot d^2, \label{i5}
\end{equation}
$$f(S)\,=\, \tfrac{2}{(\mu(S)+\mt S-1)^2 (3\mu(S)-(\mt S)^2+3\cdot\mt
S+2)^2}\,,$$
obtained in \cite{Shu96}, and
\begin{equation}
\textstyle{\sum\limits_{i=1}^r} \alpha(S_i)\,<\, 
\tfrac{2\alpha-3}{2\alpha(\alpha-1)} \cdot
d^2 - \tfrac{2\alpha-9}{2(\alpha-1)}\cdot d -
\tfrac{4\alpha}{\alpha-1}\,, 
\quad \alpha\,:=\,\max\limits_{1\le i\le r} \alpha(S_i)\,,\label{i6}
\end{equation}
where \mbox{$\alpha(\text{node})=3$}, \mbox{$\alpha(\text{cusp})=5$} and 
\mbox{$\alpha(S)\ge \frac{10}{9}(\tau'(S)+2)$} for other singularities $S$,
obtained in the Appendix to \cite{Barkats}. We like to point out that the 
coefficient of $d^2$ in (\ref{i5}) and in (\ref{i6}) depends on the ``worst'' 
singularity, hence these sufficient conditions are weakened significantly when
adding one complicated singularity. On the other hand, the new condition
(\ref{i2}) contains the contributions of the singularities in an additive form,
whence it is not so sensitive to adding an extra singularity.

\subsection*{Curves with nodes and cusps}
We pay a special
attention to the classical case of 
families of curves with $n$ nodes and $k$ cusps, for which the
criteria (\ref{i1}), (\ref{i2}) appear to be
\begin{equation}
4n+9k \,<\,d^2+6d+8\,, \quad \mbox{respectively} \quad
\tfrac{25}{2}\!\:n+18k\,<\,d^2\,, 
\label{i4}
\end{equation}
(Corollaries \ref{Corollary 3.4}, \ref{nodes and cusps}).
This is stronger than the previously known sufficient conditions
for the smoothness of 
\mbox{$\Vdi(n\!\!\;\cdot \!\!\:A_1,\ k\!\!\;\cdot\!\!\: A_2)$},
$$4n+9k<d^2\quad\mbox{(cf.\ \cite{GLS97})}\,, $$
and for the irreducibility,
$$ 225\!\:n+450\!\:k \,<\,d^2 
\quad \text{(cf.\ \cite{Shu94})}\quad \text{and} \quad
\tfrac{120}{7}\!\:n+\tfrac{200}{7}\!\:k\,<\,d^2-\tfrac{5}{7}\!\:d -
\tfrac{200}{7} \quad 
\text{(cf.\ \cite{Barkats})}\,.$$
We note also that for families of cuspidal curves our
smoothness criterion is quite close to an {\em optimal} one: 
the above inequalities
provide the smoothness and expected dimension of \mbox{$\Vdi(k\!\!\;\cdot\!\!\:
A_2)$}  for \mbox{$k\le \frac{1}{9}d^2\!\!\:+O(d)$}, 
whereas the
families of irreducible curves of degree $d$ with 
\mbox{$k=\frac{6}{49}d^2\!\!\:+O(d)$}
cusps, constructed in \cite{Shu94}, are either nonsmooth, or
have dimension greater than the expected one. That is, the coefficient 
$\frac{1}{9}$ of $d^2$ differs from an optimal one by a factor $\leq 1.1$.
 
Concerning the irreducibility it was proven in \cite{Shu94} that the 
variety \mbox{$V=V_{d}^{\text{\it irr}}(6p^2\!\cdot
  \!\!\:A_2)$} of cuspidal curves of degree $d$ with precisely
$6p^2$ cusps has at least 
two components for \mbox{$d=6p$}, showing that the coefficient $\frac{1}{18}$
of $d^2$ in (\ref{i2}) differs from an optimal one by a factor \mbox{$\leq
3$}. These examples generalize the classical example of sextic curves having 
6 cusps given by Zariski \cite{Zar}. In Proposition \ref{prop example} we 
modify the construction to obtain curves of degree $d$
slightly bigger than $6p$ having $6p^2$ cusps such that the corresponding
ESF $V$ has at least two irreducible
components but, different to Zariski's example, 
\mbox{$\pi_1(\P^2\!\setminus\!\!\:
C)=\Z/d\Z$} for each \mbox{$C\in V$}. We do not know whether $V$ is
connected.

\subsection*{Principal approach}
Looking for a sufficient smoothness and irreducibility condition, applicable to
families of curves with arbitrary singularities, we use the fact that
the smoothness and expected dimension of an equisingular 
family $V$ follow from the $h^1$-vanishing for the ideal sheaves of some
zero-dimensional subschemes of the plane (or another smooth surface) 
associated with {\it any} curve \mbox{$C\in V$} (see \cite{GrK,GrL} for a 
detailed general setting), and that the irreducibility of $V$ follows from 
the $h^1$-vanishing for the ideal sheaf of certain zero-dimensional schemes 
associated with a {\it generic} curve \mbox{$C\in V$} (such an approach 
was realized, for instance, in \cite{Shu94,Shu96,Barkats}).

Various $h^1$-vanishing criteria have been used in connection with the problems
stated. The classical idea, applied by Severi \cite{Sev},
Segre, Zariski \cite{Zar} through the later development \cite{GrK, Shu91},
is to restrict the ideal sheaf to the curve \mbox{$C\in V$} itself. 
For many cases one obtains better results when replacing $C$ by a polar curve
\cite{Shu87,GrL}, or a special auxiliary curve \cite{Shu94,Shu96}. 
A similar idea combined with Horace's method can be found in
\cite{Shu97,GLS95}.  Chiantini and Sernesi \cite{ChS} applied 
Bogomolov's theory of unstable rank two
vector bundles on surfaces for the smoothness problem of families of nodal
curves, which then was extended to curves with arbitrary singularities
\cite{GLS97}. It was Barkats \cite{Barkats} who showed how to apply 
the Castelnuovo function and Davis' Theorem \cite{Davis} for the
computation of $h^1$ in relation to the irreducibility problem. 

In the present paper we strongly exploit Barkats' observation,
combining it with other tools. Moreover, we perform our computations in a
different way to obtain stronger $h^1$-vanishing theorems (cf.\
Proposition \ref{Proposition 1} and Lemma \ref{Lemma 6.3}).
Finally, we derive
sufficient irreducibility conditions with
better asymptotic behavior (see explanation above), which involve both,
topological and analytic, singularities. A similar approach is used for
the smoothness problem completing with the result (\ref{i1}).

\subsection*{Further results and distribution of the material}
For the convenience of the reader we present the material in a self-contained 
form. In Section \ref{sec:Section 1} we introduce and set up the theory of
zero-dimensional schemes associated to singular points. In Sections 
\ref{sec:topo}
(respectively \ref{sec:anal}) we do this for topological (respectively
analytic) singularities. Section \ref{sec:hilbert} contains a proof for the
existence and irreducibility of the Hilbert scheme associated to generalized
singularity schemes, or, for clusters 
(answering a question of Kleiman and Piene). 

We compute several invariants of plane curve singularities, for instance, we
determine the degree of $C^0$-sufficiency (correcting the result in
\cite{Lichtin}), cf.\ Lemmas \ref{nu s} and \ref{Lemma 1.4}. 
In Section
\ref{sec:castelnuovo} we recall basic facts about the Castelnuovo function of a
zero-dimensional scheme in $\P^2$.  

In Sections \ref{sec:smoothness} and
\ref{sec:irred} we formulate 
the main results on the smoothness and irreducibility of equisingular families
of curves, in particular, we introduce the new invariants $\gamma(C;X)$ (cf.\
Section \ref{sec:smoothness es}). Sections
\ref{sec:Section 5} and \ref{sec:Section 6}
contain the proofs of the main results.

\subsection*{Basic definitions and notations}
Two germs \mbox{$(C,z)\subset (\P^2\!,z)$} and \mbox{$(D,w)\subset (\P^2\!,w)$}
of reduced plane curve singularities (or any of their defining power series) 
are said to be
{\em topologically equivalent} (respectively {\em analytically equivalent}, 
also called {\em contact equivalent})
if there exists a local homeomorphism (respectively analytic isomorphism)
\mbox{$(\P^2\!,z)\to (\P^2\!,w)$} mapping $(C,z)$ to $(D,w)$. The corresponding
equivalence classes are called {\em topological} (resp.\ {\em analytic})
{\em types}.  

\medskip \noindent
We recall the notion of families of plane
curves that will be used in the following.
Let $T$ be a complex space, then by a {\em family
of} (reduced, irreducible) {\em plane curves over} $T$ we mean a 
commutative diagram
$$
\renewcommand{\arraystretch}{0.8}
\arraycolsep1pt
\begin{array}{lcl}
\;\kc & \stackrel{j}{\hookrightarrow} & \;\;\P^2 \!\times T\\
\scriptstyle{\varphi} \textstyle{\searrow} &  & \swarrow \scriptstyle{pr}\\
& T &
\end{array}
$$
where $\varphi$ is a proper and flat morphism such that for all
points \mbox{$t \in T$} the fibre \mbox{$\kc_t:=\varphi^{-1}(t)$} is a 
(reduced, irreducible) plane curve, \mbox{$j : \kc
\hookrightarrow \P^2 \!\times T$} is a closed embedding and $pr$ denotes the
 natural projection.  In a similar manner, one defines {\em (flat) families of
zero-dimensional schemes} in $\P^2$ (respectively in a surface $\Sigma$). 

A {\em family with sections} is a diagram as above, together with sections 
\mbox{$\sigma_1,\dots,\sigma_r:T\to \kc$} of $\varphi$. The sections are 
called {\em trivial} if $\sigma_i$ is an isomorphism \mbox{$T \to \{z_i\}\times
  T$} for some \mbox{$z_i\in \P^2$}.
%

\medskip \noindent
To a family of reduced plane curves as above and a fibre \mbox{$C=\kc_{t_0}$}
we can associate, in a functorial way, the deformation \mbox{$\coprod_i
(\kc,z_i) \to (T,t_0)$} of the multigerm \mbox{$(C,  \text{Sing}\, C)=
\coprod_i (C,z_i)$} over the germ $(T,t_0)$. Having a family with sections
\mbox{$\sigma_1,\dots,\sigma_r$}, \mbox{$\sigma_i(t_0)=z_i$}, we obtain in the
same way a deformation  of \mbox{$\coprod_i (C,z_i)$} over $(T,t_0)$ with 
sections.  

A family \mbox{$\kc \hookrightarrow \P^2\!\times T \to T$} of reduced 
curves (with sections) is called {\em equianalytic}, respectively {\em
equisingular} (along the sections) if, for each \mbox{$t\in T$}, the induced 
deformation of the multigerm \mbox{$(\kc_t, \text{Sing}\, \kc_t)$} is 
isomorphic (isomorphic as deformation with section) to the trivial deformation,
respectively to an equisingular deformation along the trivial section
(for the equisingular case cf.\ \cite{Wahl1}). 

The Hilbert scheme of plane curves of degree $d$ together with its universal 
family is the family of all curves of degree $d$ in $\P^2$, the base space
may be identified with the linear system \mbox{$\big|
H^0\bigl(\ko_{\P^2}(d)\bigr)\big| $}. We are interested in subfamilies of
curves in $\P^2$ having fixed analytic, respectively topological types of their
singularities. 

To be specific, let \mbox{$S_1,\dots,S_r$} be fixed analytic, respectively
topological types. Denote by \mbox{$V_d(S_1,\dots,S_r)$} the space of 
reduced curves \mbox{$C\subset \P^2$} of degree $d$ having precisely $r$
singularities which are of types \mbox{$S_1,\dots,S_r$}. By \cite{GrL},
Proposition 2.1,  \mbox{$V_d(S_1,\dots,S_r)$} is a locally closed subscheme of
\mbox{$\big|H^0\bigl(\ko_{\P^2}(d)\bigr)\big| $} and represents the functor of
equianalytic, respectively equisingular families of given types
\mbox{$S_1,\dots,S_r$}. 

In the following, by abuse of notation, we write
\mbox{$C\in V_d(S_1,\dots,S_r)$} to denote either the point 
in $V_d(S_1,\dots,S_r)$ or the curve corresponding to the point, that is, the
corresponding fibre in the universal family. 

\subsection*{Acknowledgements}
We should like to thank Ragni Piene for helpful remarks concerning the
Hilbert scheme studied in Section \ref{sec:hilbert} and for the reference to the
paper \cite{NoV}. 

\smallskip
\section{Zero-dimensional schemes}
\setcounter{equation}{0}
\label{sec:Section 1}

\subsection{Geometrical meaning of zero-dimensional schemes and
  $h^1$-vanishing}\label{sec:sec1.1} \hfill\\ 
Through\-out the paper, we work with zero-dimensional schemes 
\mbox{$X=X(C)$} that are contained in a reduced plane curve \mbox{$C\subset 
\P^2$} and concentrated in finitely many points $z$.
The  corresponding ideal sheaves will be denoted
by \mbox{$\kj_{X/\P^2}\subset \ko_{\P^2}$}.
Moreover, we denote 
$$\deg X:=\textstyle{\sum\limits_z} \dim_\K
\ko_{\P^2,z}/(\kj_{X/\P^2})_z\,, \qquad \mt(X,z):=\max \left\{ \nu \in
  \Z\,\big|\, (\kj_{X/\P^2})_z\subset \fm_z^\nu \right\}\,, $$
with $\ko_{\P^2,z}$ the analytic local ring at
$z$ and \mbox{$\fm_z\subset \ko_{\P^2,z}$} the maximal ideal. 

Let $C$ be a reduced plane curve and let
\mbox{$\text{Sing}\, C =\{z_1,\dots,z_r\}$} be its singular locus.
We shall consider, among others, the following schemes $X$:
\begin{enumerate}
\itemsep3pt
\item \mbox{$\Xea(C)=\Xea(C,z_1) \cup \ldots \cup \Xea(C,z_r)$}, the
zero-dimensional  
scheme concentrated in $\text{Sing}\,C$ defined by the Tjurina ideals
$$\Iea (C,z_i)\,:=\,j(C,z_i) \,=\, \langle\, f, \tfrac{\partial f}{\partial x},
\tfrac{\partial f}{\partial y} \,\rangle \,\subset\, \ko_{\P^2,z_i}\,,$$
(where \mbox{$f(x,y)=0$} is a local equation for $(C,z_i)$). 
$\Iea(C,z_i)$ is the 
tangent space to equianalytic, i.e., analytically trivial deformations of
$(C,z_i)$.  
\item \mbox{$\Xes(C)=\Xes(C,z_1) \cup \ldots \cup \Xes(C,z_r)$}, the
zero-dimensional  
scheme defined by the equisingularity ideals 
$$\Ies (C,z_i)\,:=\,\bigl\{\, g\in \ko_{\P^2,z_i} \:\big|\: f+\varepsilon g \text{
  is equisingular over Spec}\,(\C[\varepsilon]/\varepsilon^2)\, \bigr\}\,.$$
Note that $\Xes(C)$ is contained in $\Xea(C)$ (cf.\ \cite{Wahl1}).
$\Ies(C,z_i)$ is the tangent space to equisingular deformations of $(C,z_i)$.
\item \mbox{$\Xeaf(C)=\Xeaf(C,z_1) \cup \ldots \cup \Xeaf(C,z_r)$} 
the zero-dimensional  scheme defined by the ideals
$$\Ieaf (C,z_i)\,:=\, \langle f \rangle + \mathfrak{m}_{z_i}\!\!\!\:\cdot\!\!\:
j(C,z_i) \,\subset \,j(C,z_i)\,,
$$
where \mbox{$\mathfrak{m}_{z_i}=\mathfrak{m}_{\P^2,z_i}\subset \ko_{\P^2,z_i}$}
denotes the maximal ideal.  
$\Ieaf(C,z_i)$ is the tangent space to equianalytic deformations of $(C,z_i)$
with fixed position of the singularity, i.e., equianalytic deformations along
the trivial section.
\item \mbox{$\Xesf(C)=\Xesf(C,z_1) \cup \ldots \cup \Xesf(C,z_r)$} 
the zero-dimensional  scheme defined by the ideals
$$
\Iesf (C,z_i)\,:=\,\Bigl\{ g\in \ko_{\P^2,z_i} \:\Big|\:
\renewcommand{\arraystretch}{0.6} 
\begin{array}{c} \scriptstyle{ f+\varepsilon g \text{
  is equisingular over Spec}\,(\C[\varepsilon]/\varepsilon^2) }\\ 
\text{\scriptsize along the trivial section}
\end{array}
\Bigr\}\subset  \Ies (C,z_i)\,.
$$
$\Iesf(C,z_i)$ is the tangent space to equisingular deformations of $(C,z_i)$
with fixed position of the singularity.
\item $X^s(C)=X^s(C,z_1) \cup \ldots \cup X^s(C,z_r)$
 the zero-dimensional scheme introduced in \cite{GLS96} in order to handle 
  the topological types of the singularities (cf.\ Section
  \ref{sec:topo}).
\item \mbox{$X^a(C)=X^a(C,z_1) \cup \ldots \cup X^a(C,z_r)$}
  the zero-dimensional
  scheme introduced in this paper in order to handle
  the analytic types of the singularities (cf.\ Section
  \ref{sec:anal}). In order to apply these schemes, we shall have, however, to 
  consider also (slightly) bigger schemes \mbox{$\widetilde{X}^a(C) \supset 
  X^a(C)$}.
\end{enumerate}

\smallskip
\noindent
The importance of the schemes $X(C)$ comes from the fact that the 
cohomology groups \mbox{$H^i\bigl(\kj_{X(C)/\P^2}(d)\bigr)$} have a precise 
geometric meaning for the space $V_d(S_1,\dots,S_r)$. To explain this
for $\Xeaf(C)$ and $\Xesf(C)$, consider the map
\begin{equation}
\label{Phi_d}
\Phi_d:  V_d(S_1,\dots,S_r) \lra \text{Sym}^r\P^2\,, \quad C \longmapsto
(z_1\!+\!\ldots\!+\!z_r)\,,
\end{equation} 
where $\text{Sym}^r\P^2$ is the $r$-fold symmetric product of $\P^2$ and
\mbox{$(z_1\!+\!\ldots \!+\!z_r)$} is the unordered tuple of the singularities
of $C$. 
Since any equisingular, in particular any equianalytic, deformation of a germ
admits a unique singular section (cf.\ \cite{Tei}), the universal family 
$$\ku_d(S_1,\dots,S_r)\hookrightarrow \P^2 \times V_d(S_1,\dots,S_r)
\to V_d(S_1,\dots,S_r)$$ admits, locally at $C$, $r$ singular
sections. Composing these sections with the projections to $\P^2$ gives a local
description of the map $\Phi_d$ and shows in particular that $\Phi_d$ is a well
defined morphism, even if \mbox{$V_d(S_1,\dots,S_r)$} is not reduced.

Let \mbox{$V_{d,\text{fix}} (S_1,\dots,S_r)$} denote the disjoint union of the
fibres of $\Phi_d$, together with the induced universal family on each fibre.
It follows that  \mbox{$V_{d,\text{fix}} (S_1,\dots,S_r)$} represents the
functor of equianalytic, resp.\ equisingular families of given types 
$S_1,\dots,S_r$ along trivial sections.

\medskip\noindent
In the following proposition, we write $X(C)$ instead of $\Xea (C)$, 
resp.\ $\Xeaf (C)$, resp.\ $\Xes (C)$, resp.\ $\Xesf (C)$ if the statement
holds in all four cases. Moreover, we write $V$ to denote 
\mbox{$V_{d}(S_1,\dots,S_r)$}, resp.\ \mbox{$V_{d,\text{fix}}(S_1,\dots,S_r)$}.

\begin{proposition} 
\label{H0H1}
Let \mbox{$C\subset \P^2$} be a reduced curve of degree
$d$ with precisely $r$ singularities \mbox{$z_1,\dots,z_r$} of 
analytic or topological types \mbox{$S_1,\dots,S_r$}.
\begin{enumerate}
\itemsep3pt
\item[(a)] \mbox{$H^0\bigl(\kj_{X(C)/\P^2}(d)\bigr)\big/ H^0(\ko_{\P^2})$} 
is isomorphic to the Zariski tangent space of\/ $V\!$ at $C$. 
\item[(b)] \mbox{$h^0\bigl(\kj_{X(C)/\P^2}(d)\bigr)-
h^1\bigl(\kj_{X(C)/\P^2}(d)\bigr)-1 \:\leq \: \dim (V,C) \:\leq \: 
h^0\bigl(\kj_{X(C)/\P^2}(d)\bigr)-1$}
\item[(c)]\mbox{$H^1\bigl(\kj_{X(C)/\P^2}(d)\bigr)=0$} if and only if\/ $V\!$ 
is {\em T-smooth} at $C$, i.e., 
smooth of the expected dimension \mbox{$d(d+3)/2 - \deg X(C)$}.
\item[(d)]\mbox{$H^1\bigl(\kj_{\Xea (C)/\P^2}(d)\bigr)=0$} if and only if the 
natural morphism of germs 
$$ \bigl( V_d(S_1,\dots,S_r),\,C\bigr) \lra \textstyle{\prod\limits_{i=1}^r} 
\text{Def}\,(C,z_i)$$
is smooth (in particular surjective) 
of fibre dimension \mbox{$h^0\bigl(\kj_{\Xea(C)/\P^2}(d)\bigr)-1$}.
Here \mbox{$\prod_{i=1}^r \text{Def}\,(C,z_i)$} is the cartesian product
of the base spaces of the semiuniversal deformation of the germs $(C,z_i)$.
\item[(e)] Let \mbox{$X_{\text{fix}}(C) =\Xeaf(C)$}, resp.\ $\Xesf(C)$. Then 
\mbox{$H^1\bigl(\kj_{X_{\text{fix}} (C)/\P^2}(d)\bigr)=0$} 
if and only if the morphism of germs 
\mbox{$\Phi_d: \bigl(V_d(S_1,\dots,S_r),C\bigr) \rightarrow
 \bigl(\mathrm{Sym}^r\P^2,(z_1\!+\!\ldots\!+\!z_r)\bigr) $}
is smooth of fibre dimension
  \mbox{$h^0\bigl(\kj_{X_{\text{fix}}(C)/\P^2}(d)\bigr)-1$}.
In particular, arbitrary close to $C$ there are curves in $V_d(S_1,\dots,S_r)$
having their singularities in general position in $\P^2$.
\end{enumerate}
\end{proposition}

\noindent
{\it Proof.} Note that \mbox{$H^0\bigl(\kj_{X(C)/\P^2}(d)\bigr)\big/ 
H^0(\ko_{\P^2})$} is isomorphic to 
\mbox{$H^0\bigl(\kj_{X(C)/\P^2}(d)\otimes \ko_C\bigr)$} and that
\mbox{$H^1\bigl(\kj_{X(C)/\P^2}(d)\bigr)$} is isomorphic to
\mbox{$H^1\bigl(\kj_{X(C)/\P^2}(d)\otimes \ko_C\bigr)$}. Hence the
statements (a)--(c) follow for $\Xea(C)$ and $\Xes(C)$ from \cite{GrL}, 
Theorem 3.6 (cf.\ also \cite{GrK}). The proof uses standard arguments from
deformation theory and carries over to deformations with trivial sections.
(d) was proved in \cite{GrL}, Corollary 3.9. To see (e), we apply (c) to
$X_{\text{fix}}(C)$ and notice that this implies that $\Phi_d$ has a smooth 
fibre through $C$ of the claimed dimension. Moreover,
$\kj_{X_{\text{fix}}(C)/\P^2} (d)$ is a subsheaf of $\kj_{X(C)/\P^2} (d)$,
where \mbox{$X(C)=\Xea(C)$}, resp.\ $\Xes(C)$, is of (finite) codimension $2r$.
In particular,   \mbox{$H^1\bigl(\kj_{X(C)/\P^2}(d)\bigr)=0$} and therefore,
by (c), $V_d(S_1,\dots,S_r)$ is smooth at $C$, the fibre having codimension
$2r$. It follows that $\Phi_d$ is flat with smooth fibre, hence smooth. 
\hfill $\Box$
 
\smallskip
\subsection{Zero-dimensional schemes associated to topological types of
singularities: Singularity schemes}
\label{sec:topo}

Let \mbox{$C\subset \P^2$} be a reduced plane curve of 
degree $d$ and $(C,z)$ be the
germ of $C$ at \mbox{$z\in \P^2$}, given by
\mbox{$f\in \ko_{\P^2,z}$}. 
We denote by $T(C,z)$ 
the (infinite) complete embedded resolution tree of $(C,z)$ with
vertices the points infinitely near to $z$.
We call an infinitely near point \mbox{$q\in T(C,z)$} {\it essential}, if 
it is not a node of the union of the strict
transform $f_{(q)}$ of $f$ at $q$ and the reduced exceptional divisor. 
\begin{definition}[cf.\ \cite{GLS96}]  
Let $z$ be a singular point of $C$.
We denote by $T^\ast(C,z)$ the tree spanned by $z$ and the 
essential points infinitely near to $z$.
We define $X^s(C,z)$ to be the zero-dimensional scheme given by the ideal
$$I^s(C,z)\,:=\,I^s(f)\,:=\, \bigl\{\,g \in \ko_{\P^2,z} \:\big|\:
 \mt \hat{g}_{(q)} \ge
\mt \hat{f}_{(q)},\;\; q \in T^\ast(C,z)\,\bigr\}\, \subset\, \ko_{\P^2,z},
$$
where $\hat{g}_{(q)}$ denotes the total transform of $g$ under the modification
$\pi_{(q)}$ defining $q$, and $\mt$ stands for multiplicity. We call $X^s(C,z)$
the {\em singularity scheme} of $(C,z)$.
\end{definition}

\noindent
Note that the topological type of $(C,z)$ is completely characterized 
by the partially ordered system of multiplicities 
$\mt \hat{f}_{(q)}$, \mbox{$q\in T^\ast(C,z)$}, whence 
for all elements \mbox{$g\in
I^s(C,z)$} the singularities of the germs at $z$ 
defined by $f$ and $f+tg$, $t$ generic, have the same topological type. 
Moreover, if \mbox{$g\in I^s(C,z)$} is a generic element and $(C'\!,z)$ is the
germ defined by $g$, then \mbox{$I^s(C,z)=I^s(C'\!,z)$}.
\begin{remark}
\label{remark 1.2}
We can also use the language of clusters and proximate points
(cf., e.g., \cite{Casas-Alvero}) to describe the scheme $X^s(C,z)$:
A {\em cluster} $K$ with origin at $z$ is a finite (partially ordered) set of
points $q_{i,j}$ 
infinitely near to $z$, $z$ itself included, each with assigned integral
{\em (``virtual'') 
multiplicity} $m_{i,j}$. Here, the first index $i$ refers to the {\em level} of
$q_{i,j}$, that is, the order of
the neighbourhood of $z$ which contains $q_{i,j}$.
The point \mbox{$q\in K$} is called
{\it proximate to} \mbox{$p\in K$} if it is a point in the first neighbourhood
\mbox{$E'=\pi^{-1} (p)$} of $p$, $\pi$ the blowing-up of $p$, 
or if it is a point  
infinitely near to $p$ lying on the corresponding strict transform of $E'$.
We write \mbox{$q\dashrightarrow p$}. The point \mbox{$q\in K$} is called
{\em free} if it is proximate to \mbox{$\leq 1$} point \mbox{$p\in K$}.

\smallskip \noindent
Note that for any \mbox{$q\in T(C,z)$} 
$$  \mt \hat{f}_{(q)} -\mt f_{(q)}\:= \textstyle{\sum\limits_{q
    \dashrightarrow p} }
\mt \hat{f}_{(p)}\,.  $$
Thus, it is not difficult to see that $I^s(C,z)$ is the ideal of plane curve 
germs $g$ going through the
cluster of the (partially ordered) essential points
 \mbox{$q\in T^\ast(C,z)$} with the virtual 
multiplicities $m_q:=\mt f_{(q)}$ (in the sense of \cite{Casas-Alvero},
Definition 2.3\,b). 
\end{remark}

\noindent
The degree of   $X^s(C,z)$ is
in fact an invariant of the topological type $S$ of the singularity,
namely 
$$\deg X^s(S):= \deg X^s(C,z) = \delta(C,z) + \!\textstyle{\sum\limits_{q\in
    T^\ast\!\!\:(C,z)}}  \!
m_q\,.
 $$
For this and further properties 
of $X^s(C,z)$, cf.\ \cite{GLS96} (respectively \cite{Casas-Alvero}).

\begin{definition}
Let \mbox{$(C,z)\subset (\P^2\!,z)$} be a reduced plane curve singularity
defined by \mbox{$f\in \ko_{\P^2,z}$}.
Then we define the {\em $C^0\!$-deformation-determinacy} $\nu^s(C,z)$ of 
$(C,z)$ as 
the minimum integer $\nu$ such that for any \mbox{$g\in
  \mathfrak{m}_z^{\nu+1}$} 
and all \mbox{$t\in \K$} close to $0$, the germ defined by
\mbox{$f+tg$} is topologically equivalent to $(C,z)$. 
\end{definition}

\begin{remark} 
\begin{enumerate}
\itemsep2pt
\item Recall that 
the ideal $I^s(C,z)$ defines a maximal (w.r.t.\ inclusion) linear space 
of germs $g$ such that for $t$ close to $0$ the germ 
\mbox{$f+tg$} is topologically equivalent to $(C,z)$. Hence,   
\mbox{$\nu^s(C,z)= 
\min \, \left\{ \nu \in \Z \,\big|\,  \mathfrak{m}_z^{\nu+1} \subset I^s(C,z)
\right\}$}. 
\item Let \mbox{$g\in \mathfrak{m}_z^{\nu+2}$}, \mbox{$\nu\geq
    \nu^s(C,z)$}. Then \mbox{$I^s(f+g)=I^s(f)$}. In particular, the
    singularities defined by $f$ and \mbox{$f+g$} are topologically
    equivalent.
\end{enumerate} 
\end{remark}  

\begin{lemma}
\label{nu s}
Let \mbox{$(C,z)\subset (\P^2,z)$} be a reduced plane curve 
singularity of topological type $S$ and $Q_1,\dots, Q_s$ its local 
branches. Then
$$ \nu^s(S):=\nu^s(C,z)=\min \, \Bigl\{\nu\in \Z \:\Big|\: \nu+1 \geq \max_j\,
\tfrac{2  \delta(Q_j)+ \sum_{i\neq j} (Q_i,Q_j)_z + 
\sum_{q\in T^\ast\!\cap Q_j}
  \mt Q_{j,(q)}}{\mt Q_j}\Bigr\},$$
where $(Q_i,Q_j)_z$ denotes the
intersection multiplicity of the branches $Q_i$ and $Q_j$ at $z$, and
$Q_{j,(q)}$ denotes the strict 
transform of $Q_j$ at \mbox{$q \in T^\ast\!:=T^\ast(C,z)$}. 
\end{lemma}

\noindent {\it Proof.} This follows immediately from \cite{GLS96}, Lemma 2.8.
\hfill $\Box$

\smallskip \noindent
Note that the formula for $\nu^s$ given in \cite{Lichtin} is wrong, at least in
the case of several branches, as can be seen for $A_{2k+1}$-singularities.

\bigskip \noindent
We can estimate $ \nu^s(C,z)$ in terms of
$\tau^{\text{es}}(C,z)$, the codimension 
of the $\mu$-const stratum 
in the semiuniversal deformation of $(C,z)$, respectively in terms of
$\delta(C,z)$.  Note that $\delta(C,z)$ is 
the codimension of the equiclassical stratum in the semiuniversal
deformation of $(C,z)$, whence \mbox{$\delta(C,z)\leq \tau^{\text{es}}(C,z)$}
(cf.\ \cite{DiH}).

\begin{lemma}
\label{Lemma 1.4}
\mbox{$ \nu^s(C,z) \leq  \tau^{\text{es}}(C,z)$} for any reduced plane curve
singularity $(C,z)$.
If all branches of $(C,z)$ have at least multiplicity $3$
then we have even
\mbox{$ \nu^s(C,z)  \leq  \delta(C,z)$}.
\end{lemma}

\noindent {\it Proof.} If $(C,z)$ is an $A_k$-singularity, then we have
$\tau^{\text{es}} (C,z)=\tau(C,z)=k$, and the statement is obvious. Let 
\mbox{$\mt(C,z)\geq 3$} and $Q_1,\dots,Q_s$ be the local branches of $(C,z)$.

\medskip \noindent
{\it Case 1:} $(C,z)$ is irreducible. Then, by Lemma \ref{nu s},  we have
$$ \nu^s(C,z)\:=\:\min \, \Bigl\{ \nu\in \Z \: \Big| \: \nu+1 \geq  
\tfrac{2 \delta(C,z) + \sum_{q\in T^\ast} m_q}{\mt (C,z)} \Bigr\}\,.$$
If \mbox{$\mt(C,z)=3$}, we know that \mbox{$\#\{q \in T^\ast \:\!|\!\: m_q
\!\leq\!\!\;2\}  \leq 3$}, whence  
\begin{equation}
\label{mt=3}
\textstyle{\sum\limits_{q\in T^\ast }} m_q\:\leq \:\textstyle{\sum\limits_{q\in
    T^\ast}} 
\tfrac{m_q(m_q-1)}{2} +3 \:=\:\delta(C,z)+ \mt(C,z) \,.
\end{equation}
If \mbox{$\mt(C,z)\geq 4$}, we know at least
that \mbox{$\#\{q \in T^\ast\!\: |\!\: m_q \!=\!\!\;1\} \leq \mt(C,z)$}. Thus,
$$
\textstyle{\sum\limits_{q\in T^\ast } m_q} 
\:\leq \:2 \textstyle{\sum\limits_{q\in T^\ast}} \tfrac{m_q(m_q-1)}{2} +
\mt(C,z) \:=\: 2\delta(C,z)+ \mt(C,z) \,.$$

\smallskip \noindent
{\it Case 2:} $(C,z)$ is reducible. For any $j=1,\dots,s$ we have to estimate
\begin{equation}
\label{nu est} 
\tfrac{2 \delta(Q_j)+\sum_{i\neq j} (Q_i,Q_j)_z + 
\sum_{q\in T^\ast\!\cap Q_j} \mt Q_{j,(q)}}{\mt Q_j} -1 \,.
\end{equation}
If \mbox{$\mt Q_j \geq 3$}, this does not exceed 
$$
\tfrac{
2\delta (Q_j) + \sum_{q\in T^\ast 
(Q_j)} \!\mt Q_{j,(q)}+ 2 \sum_{i\neq j} (Q_i,Q_j)_z}{\mt Q_j} -1 \:\leq \:
\delta(C,z)\,,$$
by (\ref{mt=3}) and since, as is well-known, \mbox{$\delta(C,z)= 
\sum_{i} \delta (Q_i) + \sum_{i<k} (Q_i,Q_k)_z$} (cf.\ \cite{BuG},   
Lemma 1.2.2). 
If \mbox{$\mt Q_j=2$}, (\ref{nu est}) is bounded by 
\begin{eqnarray*} 
\lefteqn{ \tfrac{1}{2} \cdot \Bigl(\textstyle{\sum\limits_{q \in T^\ast(Q_j)}
     } \mt Q_{j,(q)} 
\left(\mt Q_{j,(q)} -1\right)  + \textstyle{\sum\limits_{q\in T^\ast\!\cap
      Q_j}}  
\mt Q_{j,(q)} \left(m_q-\mt Q_{j,(q)}+1\right)\Bigr) -1 }
\hspace{1.0cm}\\
&\leq & 
\textstyle{\sum\limits_{q\in T^\ast\!\cap Q_j}} \tfrac{ m_q\cdot \mt
  Q_{j,(q)}}{2} -1 
\: \leq \: \textstyle{\sum\limits_{q\in T^\ast}} \tfrac{m_q(m_q+1)}{2} - 
\#\{q \in T^\ast\} -1 \: \leq \: \tau^{\text{es}}(C,z) 
\end{eqnarray*}
(recall that \mbox{$m_z\geq 3$} and that there are at most two points in
\mbox{$T^\ast\!\cap Q_j$} with \mbox{$m_q=1$}).
Finally, for a smooth branch \mbox{$Q_j$},   
(\ref{nu est})  can be estimated as  
$$\textstyle{\sum\limits_{q\in T^\ast\!\cap Q_j}} m_q \,-1
\:\leq \:  \textstyle{\sum\limits_{q\in T^\ast}}
 \tfrac{m_q(m_q+1)}{2} - \#\{q \in T^\ast\}-1 
\:\leq \:  \tau^{\text{es}}(C,z) 
$$
(since there is no point in \mbox{$T^\ast \!\cap Q_j$} with
\mbox{$m_q\!<2$}). 
\hfill $\Box$

\smallskip\noindent

\subsection{Hilbert schemes associated to (generalized) singularity schemes}
\label{sec:hilbert}

Let $\Sigma$ be a smooth projective surface.
The Hilbert functor $\Hilb^n_{\Sigma}$ which associates to an analytic space
$T$ the set of all (flat) families of zero-dimensional schemes in
$\Sigma$ over $T$, that is, the set of all analytic subspaces 
\mbox{$\fx \subset \Sigma \times T$}, flat over $T$ such that 
\begin{enumerate}
\item[(1)] for any \mbox{$t\in T$} the fibre $\fx_t$ of the restriction to
$\fx$ of the canonical projection \mbox{$\Sigma\times T\rightarrow T$} is
a zero-dimensional scheme of length $n$
\end{enumerate}
is well-known to be representable by a smooth connected space
of dimension $2n$, the Hilbert scheme \mbox{$\tHilb^n_\Sigma$} (cf.\
\cite{Hartshorne,Fogarty}, respectively 
the overview article \cite{Iarrobino}\,).  
That is, there is a universal family  
$$
\arraycolsep-0.1cm
\renewcommand{\arraystretch}{0.8}
\begin{array}{lcl}
\ku^n & \buildrel j\over\hookrightarrow & \;\; \Sigma \times \tHilb^n_\Sigma\\ 
\scriptstyle{\varphi} \textstyle{\searrow} & & \swarrow \scriptstyle{pr}\\
& \tHilb^n_\Sigma  & \phantom{\Big(}
\end{array}
$$
such that each element of  \mbox{$\Hilb^n_{\Sigma} (T)$}, 
$T$ a complex space, can be induced from
$\varphi$ via base change by a {\sl unique} map \mbox{$T \to
  \tHilb^n_\Sigma$}. 
Moreover, there exists a birational (``Hilbert-Chow''-) morphism 
$$ \phi: \tHilb^n_\Sigma \lra \text{Sym}^n \Sigma \,, $$ 
which can be thought of as assigning to a closed subscheme \mbox{$Z\subset
  \Sigma$} of length $n$
the $0$-cycle consisting of the points of $Z$ with multiplicities given by the
length of their local rings on $Z$ (cf.\ \cite{Fogarty}, Cor.\ 2.6).  

In \cite{Briancon}, J.\ Briancon has shown that the functor 
$\Hilb^n_{\C\{x,y\}}$, which associates to an analytic space $T$
the set of analytic subspaces \mbox{$\fx \subset \C\!\!\:^2\!\!\:\!\times T$},
 flat over $T$, satisfying (1) and
\begin{enumerate}
\item[(2)] the support of $\fx$ is contained in \mbox{$\{0\}\times T$} 
\end{enumerate}
is representable by an irreducible (but in general non-reduced)
scheme \mbox{$\tHilb^n_{\C\{x,y\}}$}. Note that 
(the reduction of) \mbox{$\tHilb^n_{\C\{x,y\}}$}
can be identified with the closed subset 
\mbox{$\phi^{-1}(nz) \subset \tHilb^n_\Sigma$}, \mbox{$z\in 
\Sigma$}. 

\begin{definition}  
Let \mbox{$\fx \hookrightarrow \Sigma\times T \to T$} be a family of
zero-dimensional schemes over a complex space $T$.
We say that the family is {\em resolvable by blowing-up sections}
if there exist pairwise disjoint sections 
\mbox{$\sigma^{(i)}_1\!,\dots,\sigma^{(i)}_{k_i}:\,
T\to Z^{(i)},$} \mbox{$i=0,\dots,N$},   
and morphisms \mbox{$\pi_i:Z^{(i+1)} \!\!\:\to Z^{(i)}$} such that
\begin{itemize}
\item \mbox{$Z^{(0)}=\Sigma\times T$}, \mbox{$\fx^{(0)}\!=\fx$},
\item \mbox{$\pi_i:Z^{(i+1)} \!\to Z^{(i)}$} is the blowup of
  $Z^{(i)}$ along the (disjoint) sections 
\mbox{$\sigma^{(i)}_1,\dots,\sigma^{(i)}_{k_i}$}
  and we denote by \mbox{$\fx^{(i+1)}$} the strict transform of 
\mbox{$\fx^{(i)}$}, \mbox{$i=0,\dots,N$}, 
\item 
for any \mbox{$0\leq i \leq N$} and any \mbox{$1\leq j \leq k_i$} the
(flat) family  
  \mbox{$(\fx^{(i)}\!\!\:\hookrightarrow \!\!\:Z^{(i)}\!\!\!\:\to\!\!\: T)$} 
is {\em equi\-mul\-tiple along the section} $\sigma^{(i)}_j$,
that is, if \mbox{$\ki_{j}^{(i)}\!\subset \ko_{Z^{(i)}}$} 
denotes the ideal of the section $\sigma^{(i)}_j$ 
then the ideal of $\fx^{(i)}$ is contained in 
$(\ki_{j}^{(i)})^m$, 
where \mbox{$m=\mt(\fx_t^{(i)}\!,\sigma_j^{(i)}(t))$} for all 
\mbox{$t\in T$},
\item \mbox{$\text{supp}( \fx^{(i)}) = \bigcup_{j=1}^{k_i}
      \sigma^{(i)}_j(T)$} and \mbox{$\text{supp}( \fx^{(N+1)}) = \emptyset$}.
\end{itemize}
\end{definition}

\noindent
\begin{remark} \label{remark 1.6}
Any (irreducible) zero-dimensional scheme $X$ supported at \mbox{$z\in
  \Sigma$} {\em defines a cluster\/} 
\mbox{$\Cl(X)$}, given by the finite set
\mbox{$ \{z\} \cup \text{supp}(X^{(1)}) \cup \ldots \cup
\text{supp}(X^{(N)})$} with assigned multiplicities
  \mbox{$m_q:=\mt(X^{(i)}\!,q)$} for 
\mbox{$q\in \text{supp}(X^{(i)})$}. 
Here, \mbox{$X^{(i+1)}\!\subset \Sigma^{(i+1)}$} is the strict transform of
  \mbox{$X^{(i)}$} under the 
  blowing-up   of \mbox{$\text{supp}(X^{(i)}) \subset \Sigma^{(i)}$},
  \mbox{$i=0,\dots,N$}, \mbox{$X^{(0)}\!:=X\subset \Sigma=:\Sigma^{(0)}$}
and \mbox{$\text{supp}(X^{(N+1)})=\emptyset$}. 
\end{remark}
                                

\noindent
Let \mbox{$(C,z)\subset(\Sigma,z)$} be a reduced plane curve 
singularity, given by \mbox{$f \in \C\{x,y\}$}, and let 
$$T^\ast\!=\{z;q_{1,1},..\!\:,
  q_{1,k_1};\dots;q_{s,1},..\!\:,q_{s,k_s}\}\subset T(C,z)$$
 be a finite subtree.
We introduce the following notations (cf.\ Remark \ref{remark 1.2}):
\begin{itemize}
\item \mbox{$K:=\Cl(C,T^\ast)$}, the cluster given by the points 
\mbox{$q\in T^\ast\!$} and the assigned virtual multiplicities 
$\underline{m}:=(\mt f_{(q)})_{q\in T^\ast\!}$;
\item \mbox{$X(C,T^\ast)$}, the zero-dimensional scheme defined by the ideal of
  plane curve germs going through the cluster \mbox{$\Cl(C,T^\ast)$};
\item \mbox{$\G:=(\Gamma_K,\underline{m})$}, the {\em cluster graph}
  associated to \mbox{$K$}.
Here $\Gamma_K$ is the (abstract)
oriented tree with coloured edges ($\lra$,$\dashrightarrow$), 
whose vertices are
in 1--1 correspondence with the points of $K$, the edges $\lra$ correspond to
pairs \mbox{$(q_{i+1,j},q_{i,k})$} with \mbox{$q_{i+1,j}$} infinitely near to
$q_{i,k}$, and the edges $\dashrightarrow$ to
pairs \mbox{$(q_{i+\ell,j},q_{i,k})$} with \mbox{$q_{i+\ell,j}$} proximate to
  $q_{i,k}$, 
\mbox{$\ell \geq 2$};
\item \mbox{$V_0$} any subset of the set of vertices $V$ of $\Gamma_K$, 
containing the root $z$ and satisfying
\begin{equation}
\label{Gamma0}
( q\in V_0\,, \; q \lra p \:\Longrightarrow\: p\in V_0 );
\end{equation}
\item \mbox{$K_0=\Cl(C,T_0^\ast)$}, \mbox{$T_0^\ast\subset T^\ast$} such that 
\mbox{$\Gamma_{K_0}=(V_0,\lra,\dashrightarrow)$}, the subgraph
of $\Gamma_K$ obtained by deleting the vertices in \mbox{$V\setminus V_0$} and
the corresponding edges;
\item \mbox{$n := \sum_{q \in V} \frac{m_q(m_q \!\!\:+\!\!\; 1)}{2}$}\,.
\end{itemize}

\noindent
Now, we define the Hilbert functor \mbox{$\Hilb^{\G,V_0}_{\kk_0}
  $} on the category of reduced
  complex spaces $T$ by associating to $T$ 
the set of all families 
\mbox{$(\fx \hookrightarrow \Sigma\times T \to T)\in \Hilb^n_{\C\{x,y\}}(T) $} 
satisfying
\begin{enumerate}
\item[($\G1$)]
there is a
finite disjoint union of irreducible reduced complex spaces \mbox{$T'\!$} and 
a finite surjective morphism \mbox{$\alpha:T'\!\to T$} 
such that the induced family \mbox{$( \alpha^\ast\fx \hookrightarrow
  \Sigma\!\!\:\times \!\!\:T'\! \to T')$} 
is resolvable by blowing-up sections \mbox{$\sigma^{(i)}_j\!:T'\to Z^{(i)}$}
(cf.\ the above definition for the notations) and, additionally, 
if $F$ is any component of the exceptional divisor of
\mbox{$Z^{(i)}\!\to \Sigma\times T'$}, then the image of $\sigma^{(i)}_j$ 
is either
contained in $F$ or it has empty intersection with $F$ (\mbox{$1\leq j\leq
  k_i$}, \mbox{$1\leq i\leq N$});
\item[($\G2$)] \mbox{$\clg(\fx_t)=\G$} for each \mbox{$t\in T$} (where
\mbox{$\clg(\fx_t)$} denotes the cluster graph defined by the cluster
$\Cl(\fx_t)$, cf.\ Remark \ref{remark 1.6}); 
\item[($\G3$)] the sections $s_j^{(i)}$ passing through the infinitely
  near points in $\Cl(\fx_t)$ corresponding to the vertices \mbox{$q_{i,j}\in
  V_0$} are trivial sections with image in \mbox{$\kk_0\!\times\!\!\; T'$}. 
\end{enumerate}

\noindent
If \mbox{$V_0$} consists precisely of the root of $\Gamma$ (i.e.,
\mbox{$\kk_0 =\{z\}$} with the assigned multiplicity $m_{z}$), we also 
write \mbox{$\Hilb^{\G}_{\C\{x,y\}}$} instead of
\mbox{$\Hilb^{\G,V_0}_{\kk_0}$}. 

\begin{remark} We use the fact, proved by A.~Nobile and O.E.~Villamayor
  \cite{NoV} (in the algebraic category), that after a finite 
base change $\alpha$ we
  always have sections, that is,
$$ \Hilb^{\G}_{\C\{x,y\}}(T)= \bigl\{ (\fx \hookrightarrow
  \Sigma\times T \to T)\in \Hilb^n_{C\{x,y\}}(T)
  \,\big|\:\clg(\fx_t)=\G \text{ for all } t\in T \bigr\}.$$
The proof of this fact can be transferred to the analytic category (cf.\ also 
\cite{Risler}). Moreover, Nobile and Villamayor show that the subfunctor 
\mbox{$\Hilb^{\G}_{\Sigma}\subset \Hilb^n_{\Sigma}$} of
families satisfying ($\G1$) and
($\G2$) (defined on the category of reduced algebraic schemes) 
is representable.
\end{remark}

\begin{proposition}
\label{locally closed} 
The functor \mbox{$\Hilb^{\G,V_0}_{\kk_0}$} 
(defined on the category of reduced complex spaces) 
is representable
by a locally closed subspace 
\mbox{$\tHil^{\G,V_0}_{\kk_0}\subset
\tHilb^n_{\C\{x,y\}}$}.
In particular, the functor $\Hilb^{\G}_{\C\{x,y\}}$ is representable by a 
locally closed subspace 
\mbox{$\tHilb^{\G}_{\C\{x,y\}}\subset \tHilb^n_{\C\{x,y\}}$}. 
\end{proposition}
 
\begin{proof} 
We proceed by induction on $n$. For \mbox{$n=1$} there is nothing to show
(\mbox{$\tHilb^1_{\C\{x,y\}}$} is just one point). Let \mbox{$n\geq 2$} and 
\mbox{$m:=m_{z}$}.

\medskip\noindent
{\it Step 1.} The subfunctor $\kh^{n,m}_{\C\{x,y\}}$ of
$\Hilb^{n}_{\C\{x,y\}}$ given by
$$ \kh^{n,m}_{\C\{x,y\}} (T) = \bigl\{(\fx \hookrightarrow \Sigma\times T \to
T)\in \Hilb^{n}_{\C\{x,y\}} (T) 
\:\big|\: \mt(\fx_t)=m \text{ for any } t\in T\/\bigr\} $$
is representable by a locally closed subspace
\mbox{$H^{n,m}_{\C\{x,y\}} \subset \tHilb^n_{\C\{x,y\}}$}.

\medskip \noindent
This can be seen as follows: Consider the description of \mbox{$
\tHilb^n_{\C\{x,y\}}$} as an algebraic subset of the Grassmannian of codim
$n$ vector spaces of \mbox{$\C\{x,y\}/\mathfrak{m}^n$} given by J.~Briancon.
In the local coordinates $\lambda_{\alpha,\beta,i,j}$ 
associated to given stairs 
(cf.\ \cite{Briancon}, II\,2.1) the subspace
\mbox{$H^{n,m}_{\C\{x,y\}}\subset\tHilb^n_{\C\{x,y\}}$}
is defined by the vanishing of all $\lambda_{\alpha,\beta,i,j}$ with \mbox{$i+
j<m$} and the condition that not all $\lambda_{\alpha,\beta,i,j}$, \mbox{$i+j=m
$} vanish.

\medskip\noindent
{\it Step 2.} We introduce the following notations:
\begin{enumerate}
\item[(a)]
\mbox{$\fu \hookrightarrow \Sigma\times
  H^{n,m}_{\C\{x,y\}}\to H^{n,m}_{\C\{x,y\}}$} denotes the universal family,
\mbox{$\pi:\Sigma'\to \Sigma$}, respectively 
\mbox{$\pi:\Sigma'\!\times H^{n,m}_{\C\{x,y\}}\to 
\Sigma\times H^{n,m}_{\C\{x,y\}}$} the blowing-up of \mbox{$z\in
\Sigma$} (respectively of the trivial section \mbox{$t \mapsto (z,t)$}
in \mbox{$\Sigma\times H^{n,m}_{\C\{x,y\}}$}), and \mbox{$E'$}
the exceptional divisor of $\pi$.  
\item[(b)] \mbox{$\G^{(i)}\!:=(\Gamma^{(i)},\underline{m})$} \, 
(\mbox{$\Gamma^{(i)}$} an oriented coloured tree with set of vertices
$V^{(i)}$), \mbox{$i=1,\dots,k_1$},
denote the cluster graphs obtained from $\G$ by removing the root $z$ (cf.\
  Figure \ref{fig trees}). Set   
\begin{figure}
\centerline{\psfig{figure=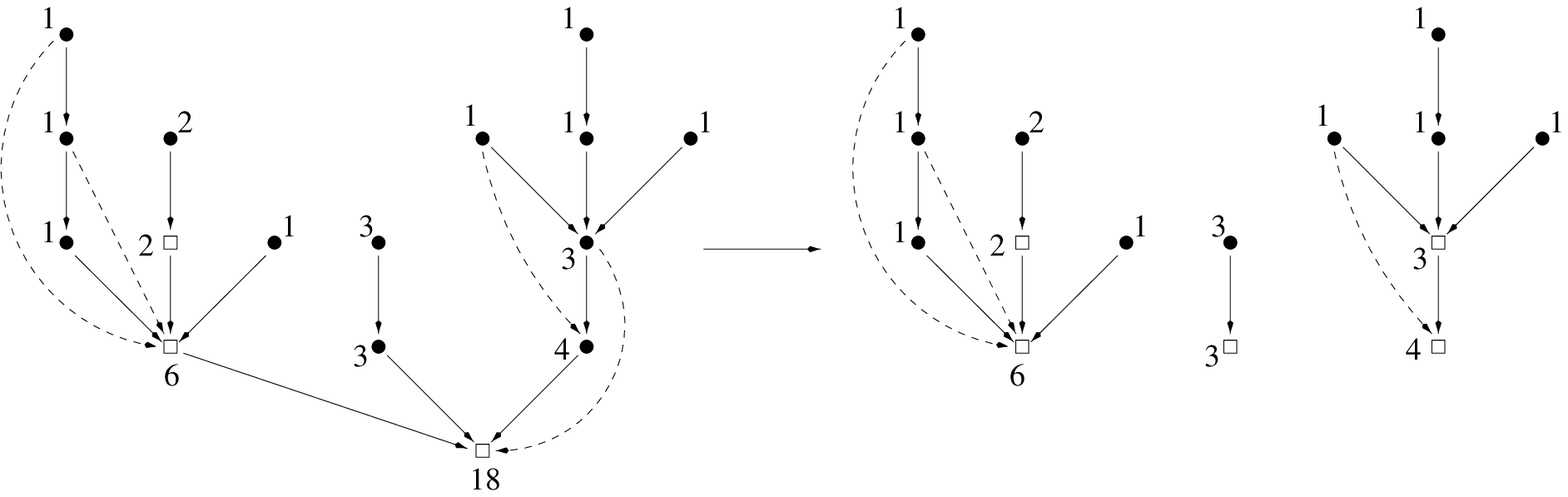,height=4.2cm,width=12cm}}
\caption{\label{fig trees} A cluster graph $\G$ (with subset 
\mbox{$V_0\subset V$}, marked by $\Box$) and the 
cluster graphs $\G^{(i)}$ (with subsets \mbox{$V^{(i)}_0\!\subset V^{(i)}$}), 
\mbox{$i=1,2,3$}.}
\end{figure} 
$$ n_i:= \sum_{q\in V^{(i)}} \frac{m_q(m_q\!+\!\!\:1)}{2}\,,\quad 
\widetilde{n}:=n_1+\ldots + n_{k_1} = n-\frac{m(m\!+\!\!\:1)}{2}\,.$$
\end{enumerate}

\noindent
Without restriction, we can assume that the roots of
\mbox{$\Gamma^{(1)}\!,\dots,\Gamma^{(s)}$}, \mbox{$0\leq s \leq k_1$}, 
are vertices in $V_0$ (corresponding precisely to the infinitely near points 
$q_{1,1},\dots,q_{1,s}$ of level 1 in $\kk_0$), 
while the roots of
$\Gamma^{(i)}$, \mbox{$i>s$}, are not in $V_0$. 
We introduce the subsets 
$$V_0^{(i)}:= \bigl(V_0\cap V^{(i)}\bigr)  \cup 
\underbrace{\{\text{root of $\Gamma^{(i)}$}\} \cup
\bigl(\left \{p \mid p \dashrightarrow
    z\right\} \cap \, V^{(i)}\bigr)}_{ 
\renewcommand{\arraystretch}{0.7}
\begin{array}{c}
    \scriptscriptstyle{||}\\ 
\textstyle{\{p_{\ell_i,i} \to \dots \to p_{1,i} \}}
\end{array}  }
\subset V^{(i)} \,
,\quad i=1,\dots,k_1,$$ 
which (clearly) satisfy the property (\ref{Gamma0}) and which correspond to 
clusters $\kk^{(i)}_0$  with origin  
\mbox{$p_{1,i}\in E'\subset \Sigma'$}, \mbox{$i=1,\dots,k_1$}, given by
\begin{itemize}
\item those points in \mbox{$\kk_0$} which are infinitely
  near to  \mbox{$p_{1,i}=q_{1,i}$}, \mbox{$i=1,\dots,s$},    
\item  the intersection points $p_{j,i}$, \mbox{$j=2,\dots,\ell_i$}, 
of the strict transform of $E'$ 
with the exceptional divisor of \mbox{$\pi_{j,i}:\Sigma_{j,i}\to
  \Sigma_{j-1,i}$}, the blowup of $p_{j-1,i}$ in
$\Sigma_{j-1,i}$ (where \mbox{$\Sigma_{1,i}=\Sigma'$}), \mbox{$i=1,\dots,k_1$}.
\end{itemize}
Note that the points \mbox{$p_{1,1},\dots,p_{1,s}$} are already fixed by
$\kk_0$, while \mbox{$p_{1,s+1},\dots,p_{1,k_1}$} can be chosen arbitrarily in
$E'$, such that all the $p_{1,i}$ are pairwise distinct.

\medskip \noindent
{\it Step 3.} Let \mbox{$t\in H$} be such that \mbox{$\clg(\fu_t)=\G$} 
and such that the infinitely near points corresponding to the vertices in 
$V_0$ are in the prescribed position given by $\kk_0$. We
show that there exists a cartesian diagram of germs
$$
\SelectTips{cm}{12}
\xymatrix@C=5pt@R=8pt{
 \bigl(H^{n,m}_{\C\{x,y\}},t\bigr)\phantom{\Big|}\\
 (H,t)\phantom{\Big|}\ar[r]^-{\psi}_-{(i)} \ar@{^{(}->}[u]_{\text{closed}}&
\,\bigl(\tHilb^{\widetilde{n}}_{\Sigma'},\psi(t)\bigr)
\,\ar[r]^{\zeta}_{\cong} 
&
\,\textstyle{\prod\limits_{i=1}^{k_1}}
 \bigl(\tHilb^{n_i}_{\Sigma'}, \psi(t)_i \bigr)\\
 & & \hspace{-1.5cm}\textstyle{\prod\limits_{i=1}^{s}}
 \bigl(\tHilb^{n_i}_{\C\{x,y\}}, \psi(t)_i \bigr) \!\!\:\times\!\!\!\:
 \textstyle{\prod\limits_{i=s+1}^{k_1}} 
 \bigl(E'\!\times\!\!\:\tHilb^{n_i}_{\C\{x,y\}}, (p_{1,i},\psi(t)_i) \bigr)
\phantom{\Big|}\ar@{^{(}->}[u]^{(ii)\,}_{\text{closed}}\\
 \bigl(\tHil^{\G,V_0}_{\kk_0},t\bigr)\phantom{\Big|} 
\ar@{^{(}->}[uu] \ar[r]_-{\zeta\circ\psi}
&\,\textstyle{\prod\limits_{i=1}^{s}}\phantom{\Cl} & 
\hspace{-1.6cm}
 \bigl(\tHil^{\G^{(i)},V^{(i)}_0}_{\kk^{(i)}_0}\!, \psi(t)_i \bigr)
 \!\!\:\times\!\!\!\: 
 \textstyle{\prod\limits_{i=s+1}^{k_1}} 
 \bigl(E'\!\times\!\!\:\tHil^{\G^{(i)},V^{(i)}_0}_{\kk^{(i)}_0}\!,
 (p_{1,i},\psi(t)_i) \bigr) 
\ar@{^{(}->}[u]^{(iii)\,}_{\text{closed}}
} 
$$
obviously implying the statement of Proposition \ref{locally closed}.

\medskip
\noindent
\mbox{$(i)$}\, We consider the strict transform 
\mbox{$\varphi:\widetilde{\fu}\hookrightarrow \Sigma'\!\times\!\!\:
  H^{n,m}_{\C\{x,y\}}\! \to H^{n,m}_{\C\{x,y\}}$} of the universal family,  
given by the ideal (sheaf) $\kj_{\widetilde{\fu}}$ associated to
\mbox{$U\mapsto \kj_{\widetilde{\fu}}(U):= \bigl\{\;\!\widehat{g}
\;\big|\; g\in \!\!\:\kj_{\fu}(U) \bigr\}:(\ki_{E'}(U))^m$}. Here,
$\widehat{g}$ denotes the total transform of $g$ under $\pi$, and $\ki_{E'}$ 
the ideal of the exceptional divisor in \mbox{$\Sigma'\!\times\!\!\:
  H^{n,m}_{\C\{x,y\}}$}. 

By semicontinuity of the fibre dimension of the finite morphism
$\varphi$, it follows that there is a locally closed subset \mbox{$H\subset
 H^{n,m}_{\C\{x,y\}}$} such that for any \mbox{$t\in H$} we have
\mbox{$ \dim_{\C} (\widetilde{\fu}_{\!\;t})= \widetilde{n}$}.
In particular, the restriction 
of $\varphi$ to the preimage of $H$ defines a flat morphism, whence, by the 
universal property of $\tHilb^{\widetilde{n}}_{\Sigma'}$ there exists a
morphism
\mbox{$\psi: H\to \mathrm{Hilb}^{\widetilde{n}}_{\Sigma'}\,. $}

\medskip
\noindent
$(ii)$\,
There is an isomorphism of germs
\mbox{$\zeta:\,\bigl(\mathrm{Hilb}^{\widetilde{n}}_{\Sigma'}, \psi(t)\bigr) 
\stackrel{\cong}{\lra} 
\textstyle{\prod_{i=1}^{k_1}}
 \bigl(\mathrm{Hilb}^{n_i}_{\Sigma'}, \psi(t)_i) $}
(cf., e.g., \cite{Iar}), and we can consider the (Hilbert-Chow) 
morphism of germs
$$ \phi=(\phi_1,\dots,\phi_{k_1}): \: \textstyle{\prod\limits_{i=1}^{k_1}}
 \bigl(\mathrm{Hilb}^{n_i}_{\Sigma'}, \psi(t)_i) \lra
\textstyle{\prod\limits_{i=1}^{k_1}}
 \bigl(\mathrm{Sym}^{n_i} \Sigma', n_i\!\cdot \!\!\:p_{1,i}\bigr)\,.
$$
The preimages under $\phi_i$ of the (germs 
at \mbox{$n_i\!\cdot \!\!\:p_{1,i}$} of the) 
locally closed subsets 
$$ \Delta^{(i)} := \Biggl\{
\begin{array}{ll}
\left\{n_i\!\cdot \!\!\:p_{1,i}\right\}  & \text{ if }\, 1\leq i \leq s \\
\left\{n_i\!\cdot \!\!\:w\mid w\in E' \right\}  & \text{ if }\, s < i \leq
k_1 
\end{array}
$$ 
are (locally) isomorphic to  
$\tHilb^{n_i}_{\C\{x,y\}}$ (if \mbox{$i\leq s$}),
respectively to \mbox{$E'\!\times\tHilb^{n_i}_{\C\{x,y\}}$} (if
\mbox{$i>s$}). 

\medskip
\noindent
$(iii)$\,
Finally, locally at $t$,
\mbox{$\tHil^{\G,V_0}_{\kk_0}$} is the preimage under
\mbox{$\zeta\circ \psi$} of 
$$ \textstyle{\prod\limits_{i=1}^s}
 \tHil^{\G^{(i)},V^{(i)}_0}_{\kk^{(i)}_0} \times \!
\textstyle{\prod\limits_{i=s+1}^{k_1}}
\! \Bigl(E'\!\times\!\!\; 
\tHil^{\G^{(i)},V^{(i)}_0}_{\kk^{(i)}_0}\Bigr)
\:\subset \; \textstyle{\prod\limits_{i=1}^s}
 \tHilb^{n_i}_{\C\{x,y\}} \times \!
\textstyle{\prod\limits_{i=s+1}^{k_1}}
\! \left( E'\!\times\!\!\; \tHilb^{n_i}_{\C\{x,y\}}\right)
$$
which, by the induction hypothesis, is a locally closed subset. 
\end{proof}

\begin{proposition}
\label{irreducible}
The Hilbert scheme \mbox{$\tHil^{\G,V_0}_{\kk_0}$} is irreducible and has 
 dimension $M$ equal to the number of free points in \mbox{$K\setminus K_0$}. 
In particular, \mbox{$\tHilb^{\G}_{\C\{x,y\}}$} is irreducible of
 dimension equal to the number of free points in \mbox{$K\setminus \{z\}$}.  
\end{proposition}

\begin{proof}  
Again, we proceed by induction on $n$. With the notations introduced in the
proof of Proposition \ref{locally closed}, we can assume that the first $\ell$
triples  
$$\bigl(\G^{(i)}, V_0\cap V^{(i)}, V_0^{(i)}\bigr)\,, \quad  
i=1,\dots,\ell\,,$$
are pairwise different and occur precisely $\nu_i$-times among all such triples
(in particular, \mbox{$\nu_1+\ldots+\nu_\ell=k_1$}). Recall that we assumed
\mbox{$V_0\cap V^{(i)} \neq\emptyset$} precisely 
for \mbox{$i=1,\dots,s\leq \ell$}. (Note that \mbox{$\nu_i=1$}
if \mbox{$V_0\cap V^{(i)} \neq\emptyset$}).

For any \mbox{$i=1,\dots,\ell$}, let $\widetilde{\fx}^{(i)}$
be the union of those connected components of the strict transform
\mbox{$\varphi:\widetilde{\fx}\hookrightarrow \Sigma'\!\times
  \!\!\:\tHil^{\G,V_0}_{\kk_0}\!\to 
\tHil^{\G,V_0}_{\kk_0}$} 
of the universal family
which satisfy 
\begin{itemize}
\item \mbox{$\clg(\widetilde{\fx}^{(i)}_t\!,x)=\G^{(i)}$},
\item the infinitely near points of \mbox{$\Cl(\widetilde{\fx}^{(i)}_t\!,x)$}
 corresponding to the vertices in
  \mbox{$V_0\cap V^{(i)}$} are in the prescribed position given by $\kk_0$,
\item the infinitely near points of \mbox{$\Cl(\widetilde{\fx}^{(i)}_t\!,x)$}
 corresponding to the vertices in \mbox{$V_0^{(i)}$} are on $E'$ (respectively
 on its strict transform)
\end{itemize}
for all \mbox{$x\in \text{supp}(\widetilde{\fx}^{(i)}_t)$}, 
\mbox{$t\in \tHil^{\G,V_0}_{\kk_0}$}.
In particular, 
\mbox{$\widetilde{\fx} =\widetilde{\fx}^{(1)}\cup \ldots \cup
\widetilde{\fx}^{(\ell)}$} and the fibres of the restriction of $\varphi$,\, 
\mbox{$\varphi_i:\widetilde{\fx}^{(i)}\to \tHil^{\G,V_0}_{\kk_0}$}
have constant (vector space) dimension (\mbox{$=\nu_in_i$}).
Hence the $\varphi_i$ are flat and, by the universal property of
\mbox{$\tHilb^{\nu_in_i}_{\Sigma'}$}, we obtain morphisms
$$ \tHil^{\G,V_0}_{\kk_0} \stackrel{\rho_i}{\lra}
\tHilb^{\nu_i n_i}_{\Sigma'} 
\stackrel{\phi_i}{\lra} \text{Sym}^{\nu_in_i} \Sigma'\,, \quad
i=1,\dots,\ell\,.$$ 
We complete the proof by showing that the composed morphism
$$\phi\circ \rho:= (\phi_1\circ \rho_1,\dots,\phi_\ell\circ \rho_\ell):\;
\tHil^{\G,V_0}_{\kk_0} \lra \text{Sym}^{\nu_1 n_1} \Sigma'
\times \ldots \times \text{Sym}^{\nu_\ell n_\ell} \Sigma'$$  
is dominant with irreducible and equidimensional fibres on the irreducible set 
\mbox{$\Delta_1\times \ldots\times \Delta_\ell$}.
Here, \mbox{$\Delta_i=\big\{n_i\!\cdot \!\!\:q_{1,i}\big\}$} 
if \mbox{$1\leq i\leq s$}  (\mbox{$q_{1,i}$} being the
infinitely near point in $\kk_0$ corresponding to the root of
\mbox{$\Gamma^{(i)}$}), 
and \mbox{$\Delta_i=\big\{\sum_{j=1}^{\nu_i} n_i\!\cdot \!\!\:w_{i,j}\,\big|\,
w_{i,j}\in E'\big\}$} if \mbox{$s< i\leq \ell$}.

Let $(w_{i,j})_{i,j}$ be any $k_1$-tuple 
of {\em pairwise different} points \mbox{$w_{i,j}\in E'$}, 
\mbox{$w_{i,1}=q_{1,i}$} if  \mbox{$1\leq i\leq s$},
(\mbox{$j=1,\dots,\nu_i$}, \mbox{$i=1,\dots,\ell$}).
Then there is a curve germ $(C(\underline{w}),z)$, topologically
equivalent to $(C,z)$, having tangent directions $w_{i,j}$. Moreover, we can
choose $C(\underline{w})$ such that the local branches of $C$ and
$C(\underline{w})$  with tangent direction \mbox{$q_{1,i}$},
\mbox{$i=1,\dots,s$}, 
coincide. By chosing the subtree \mbox{$T^\ast(\underline{w})\subset
  T(C(\underline{w}),z)$} corresponding to \mbox{$T^\ast\!\subset
  T(C,z)$}, we obtain a zero-dimensional scheme
\mbox{$X(\underline{w})=X(C(\underline{w}),T^\ast(\underline{w}))$} with
associated cluster graph $\G$. By construction, \mbox{$X(\underline{w})$}
corresponds to a point in the fibre 
\mbox{$(\phi\circ \rho)^{-1} \bigl( \sum_{j=1}^{\nu_1}
    n_1 w_{1,j}, \dots, 
\sum_{j=1}^{\nu_\ell} n_\ell w_{\ell,j} \bigr)$}.
On the other hand, any point in the image is of this form and 
$$ (\phi\circ \rho)^{-1} \Bigl(\textstyle{\sum\limits_{j=1}^{\nu_1} n_1
  w_{1,j}, \dots, 
\sum\limits_{j=1}^{\nu_\ell} n_\ell w_{\ell,j}\Bigr)  \,\cong \,
\prod\limits_{j=1}^{\nu_1} \tHil^{\G^{(1)},V^{(1)}_0}_{\kk^{(1)}_0}
\times \ldots \times  \prod\limits_{j=1}^{\nu_\ell}
\tHil^{\G^{(\ell)},V^{(\ell)}_0}_{\kk^{(\ell)}_0}}\,.$$
Hence, by the induction hypothesis, the fibres are irreducible and 
  equidimensional.  

In the same manner, 
the dimension statement follows from the induction hypothesis, since the
dimension of the image of \mbox{$\phi\circ \rho$} equals the number of free
points of level 1 in \mbox{$K\setminus K_0$}. 
\end{proof}

\begin{remdef} \label{1.10}
Let \mbox{$(C,z)\subset (\Sigma,z)$} be a reduced plane curve singularity. 
Then, by the above, the cluster graph $\G$ defined by the cluster
\mbox{$\Cl(C,T^\ast(C,z))$} is an invariant of
the topological type $S$ of the singularity. Hence, we can introduce
\label{kx0t} 
$$ \kxf(S) := \tHilb^{\G}_{\C\{x,y\}}\,. $$
Notice that the universal family \mbox{$\ku_d(S)\hookrightarrow
\P^2\!\times \!\!\:V_d(S) \rightarrow V_d(S)$}  of reduced plane curves of
degree $d$ 
having a singularity of (topological) type $S$ along the section 
\mbox{$\Phi_d: V_d(S)\to \P^2$} 
as its only singularity 
defines a family \mbox{$\varphi: \fx^s \hookrightarrow \P^2\!\times\!\!\:
V_d(S) \rightarrow V_d(S)$} of singularity schemes (supported along $\Phi_d$).
There exists 
an affine subset \mbox{$\A^2\!\subset \P^2$} such that the complementary line
\mbox{$L_\infty$} satisfies
$$V:= V_d(S)\setminus \Phi_d^{-1}(L_\infty)
\underset{\text{dense}}{\hookrightarrow}  V_d(S)\,.$$
We consider the induced family \mbox{$\fx^s\!\!\: \hookrightarrow \A^2\!\times
V \rightarrow V$}. Applying the translation 
$$ \A^2\!\times V \lra \A^2\!\times V: \; (\underline{x};\,C) \longmapsto
\bigl(\underline{x}\!\!\:-\!\!\:\Phi_{d}(C);\,C) $$
leads to a family over $V$ of zero-dimensional schemes in $\A^2$, supported
along the trivial section. It follows that there exists 
a morphism
\begin{equation}
\label{Psi_d 2}
 \Psi_d: V  \lra \kx(S):=\P^2\!\times \kxf(S) \,,
\end{equation}
assigning to a curve \mbox{$C\in V$} with singularity at $w$ the tuple 
\mbox{$\bigl( \Phi_d(C),\tau_{w0}(X^s(C,w))\bigr)$},
where $\tau_{w0}$ denotes the translation mapping $w$ to $0$.
\end{remdef}

\subsection{Zero-dimensional schemes associated to analytic types of singular
points}
\label{sec:anal}

Even if throughout the paper we work with plane curves, we 
should like to introduce
the analogue to the schemes $X^s$ for analytic types in the more general 
context of hypersurfaces \mbox{$F\subset \P^n$} with isolated singularities.

Let \mbox{$f\in \ko_{\P^n\!,z}$} define an isolated singularity.  
We consider zero-dimensional ideals $I(g) \subset \ko_{\P^n\!,z}$ defined 
for every \mbox{$g \in \ko_{\P^n,w}$}  analytically (or contact) equivalent to
$f$, that is,  of the form \mbox{$g=(u\!\!\:\cdot\!\!\: f) \circ \psi$}
 with \mbox{$\psi:(\P^n\!,w) \rightarrow (\P^n\!,z)$} 
a local analytic isomorphism and \mbox{$u\in \ko_{\P^n\!,z}$} a unit, 
such that the following four conditions hold:
\begin{enumerate}
\itemsep1pt
\item[(a)] \mbox{$g\in I(g)$},
\item[(b)] a generic element \mbox{$h\in I(g)$} is contact equivalent to $g$
and satisfies \mbox{$I(h)=I(g)$},
\item[(c)] for $\psi$ and $u$ as above we have  
\mbox{$I(\psi^\ast(u\cdot f))=\psi^\ast  I(f)$}. 
\item[(d)] there exists an \mbox{$m>1$} such that $I(g)$ is determined by the
$m$-jet of $g$.
\end{enumerate}
Note that (c) implies that this definition is independent of the
choice of the generator $g$ of the ideal $\langle g \rangle$ and that
the isomorphism class of $I(g)$ is an invariant of the analytic type of $g$.
If the germ \mbox{$(F,z)\subset (\P^n\!,z)$} is given by $f$, we set
\mbox{$I(F,z):=I(f)$} and \mbox{$X(F,z):=V(I(F,z))\subset \P^n$}.

\begin{definition}
Let \mbox{$(F,z)\subset (\P^n\!,z)$} be a hypersurface germ with isolated 
singularity given by
\mbox{$f\in \ko_{\P^n\!,z}$}. If a collection of ideals $I(g)$, $g$ contact
equivalent to $f$, satisfies
(a)--(d) and has the maximal possible size, i.e., minimal colength in 
$\ko_{\P^n,w}$, we denote $I(g)$ by $I^a(g)$. We set
$$ I^a(F,z):= I^a(f)\,, \quad X^a(F,z)=V( I^a(F,z)) \subset \P^n \,.$$
Since the degree of the zero-dimensional scheme $X^a(F,z)$ is invariant under
local analytic isomorphisms we can introduce \mbox{$\deg X^a(S):=\deg
X^a(F,z)$}, where $S$ is the analytic type of $(F,z)$.
Moreover,
since \mbox{$X^a(F,z)=X^a(f)$} is zero-dimensional, we can define 
$$ \nu^a(F,z):=\nu^a(f):= \min \, \bigl\{ \nu \in \Z
\,\big|\,  \mathfrak{m}_z^{\nu+1} \subset
I^a(f) \bigr\}\,.$$
$\nu^a(F,z)$ is called the {\em (analytic) de\-for\-mation-de\-ter\-mi\-nacy} 
of
$(F,z)$.  Note that $\nu^a$ does only depend on the analytic type $S$ of the
singularity $(F,z)$. Hence, we may introduce $\nu^a(S):=\nu^a(F,z)$.
\end{definition}

\noindent
Recall that the analytic type of an isolated hypersurface
singularity \mbox{$(F,z)\subset (\P^n\!,z)$} with Milnor
number \mbox{$\mu=\mu(F,z)$} is already determined by its
\mbox{$(\mu\!+\!1)$}-jet. Hence, by the maximality of $I^a(f)$,
\mbox{$\nu^a(F,z)\leq \mu(F,z)+1$}. 
We shall show that even \mbox{$\nu^a(F,z)\leq \tau(F,z)$}, where $\tau(F,z)$
denotes the Tjurina number of $(F,z)$. 
\begin{remark}
\label{M irred}
Let $S$ be an analytic type,
\mbox{$0=(0\!\!\::\!\!\:\ldots \!\!\::\!\!\:0\!\!\::\!\!\:1)\in \P^n$} and
\mbox{$f\in \ko_{\P^n\!,0}$} 
define a singularity of type $S$. 
Consider a collection of 
ideals $I(g)$, $g$ contact equivalent to $f$, satisfying (a)--(d).
The set of all
zero-dimensional schemes \mbox{$X(F,0)\subset \P^n$}, 
$(F,0)$ being of type $S$, 
coincides with the set of all $X(g)$, \mbox{$g\in \ko_{\P^n,0}$} contact
equivalent to $f$, which,
by condition (c), can be identified with the orbit of $I(f)$ 
(mod $\mathfrak{m}_0^{\nu+1}$) under 
the action of the (irreducible) algebraic group 
$$G\:=\: \text{Diff}\,\text{(mod}\,\mathfrak{m}_{0}^{\nu+1}).$$
Here Diff denotes the group of local analytic isomorphisms
\mbox{$(\P^n,0) \rightarrow (\P^n,0)$} and \mbox{$\nu \geq \nu^a(f)$}. 
\end{remark}

\begin{definition}
Let $\kxf(S)$ denote the orbit of $I(f)$ 
(mod $\mathfrak{m}_0^{\nu+1}$) under the action of $G$.  
Let $V_d$ be (the base space of) a family of reduced 
hypersurfaces $F$ of degree $d$  having an isolated singularity of type $S$ 
along the section \mbox{$z=z(F)$}. As in Remark \ref{1.10}, there exists
\mbox{$\A^n\!\subset \P^n$} and a dense 
subset \mbox{$V\subset V_d$} such that the support of $X(F)$, \mbox{$F\in V$},
is contained in $\A^n$. In particular, we can define a morphism
\begin{equation}
\label{Psi_d 1}
 V_d \underset{\text{dense}}{\supset} 
V \stackrel{\Psi_d}{\lra}  \kx(S):=\P^n\!\times \kxf(S)\,,\quad
F  \longmapsto  \left(z,\tau_{z0}(X^a(F,z))\right)  \,,
\end{equation}
where $\tau_{z0}$ denotes the translation mapping $z$ to $0$. Note that 
$\kx(S)$ is irreducible by Remark \ref{M irred}.
\end{definition}

\smallskip\noindent
In general, the schemes $X^a(F,z)$ are difficult to handle, since there is no 
concrete description of $I^a(F,z)$, which would be needed, 
e.g., to determine the degree of $X^a(F,z)$. 
Of course, there are special cases, where we can describe $I^a(F,z)$ 
explicitely. For instance, for a simple plane curve singularity
$(C,z)$, where we have just \mbox{$X^a(C,z)=X^s(C,z)$}. 

\smallskip \noindent
To be able to estimate
$\deg X^a(S)$ for arbitrary singularities we shall introduce
ideals $I(g)$ satisfying the properties (a)--(d), but not necessarily being
of maximal size.

Note that necessarily \mbox{$I(g) \subset \Ieaf(g)=\langle g \rangle + 
\mathfrak{m}_z \!\cdot \!\!\;j(g)$}, since for \mbox{$h \in I(g)$} the
deformation 
\mbox{$g+th$} is equianalytic with fixed position of the singularity, in
particular, the tangent vector $h$ to this deformation is an element of
$\Ieaf(g)$.

\begin{definition}
Let \mbox{$f\in \ko_{\P^n\!,z}$} be an isolated singularity and 
let $j(f)$ denote the Tjurina ideal, i.e., the ideal generated by $f$ and its
partial derivatives. 
We introduce
$$
\widetilde{I}^a(f)\::= \:\bigl\{ \, g\in \ko_{\P^n\!,z} \;\big| \; 
                       \: j(g) \subset j(f)\,\bigr\} \:\subset\: j(f)\,.
$$
If \mbox{$\underline{x}=(x_1,\dots,x_n)$}
are local coordinates at $z$ 
and if \mbox{$f\in \K\{\underline{x}\}$} then 
\begin{equation}
\label{tilde I}
\widetilde{I}^a(f)\:= \:
\left\{ \alpha_0 f + \textstyle{\sum\limits_{i=1}^n} \alpha_i 
\tfrac{\partial f}{
\partial x_i} \;\,\Big| \,\; \begin{array}{l} \alpha_0,\,
  \alpha_1,\dots, \alpha_n \in \K\{\underline{x}\}\,,\\
(\alpha_1,\dots,\alpha_n)\cdot D^2\! f(\underline{x}) \equiv 
\underline{0} \text{ mod } j(f) 
\end{array}
\right\} 
\end{equation}
where $D^2 \!f(\underline{x})$ denotes the Hessian matrix. 
\end{definition}

\noindent
Clearly, $\widetilde{I}^a(f)$ is an ideal containing $f$ and it is already 
determined by the \mbox{$(\mu\! +\!1)$}-jet of $f$.
We shall show that the collection of ideals
$\widetilde{I}^a(f)$ satisfies
also the conditions $(b)$ and $(c)$. The description (\ref{tilde I}) of 
$\widetilde{I}^a(f)$ provides an algorithm, using standard bases, to compute
$\widetilde{I}^a(f)$, which has been implemented in {\sc Singular} \cite{GPS},
cf.\ \cite{Lo1}.

\begin{lemma and definition}
Let \mbox{$z,w\in \P^n$} be arbitrary points. Moreover, let
\mbox{$f\in \ko_{\P^n\!,z}$} be an isolated singularity,
\mbox{$\psi: (\P^n\!,w) \to (\P^n\!,z)$} the germ of an analytic isomorphism
and $ u\in \ko_{\P^n\!,z} $ a unit. Then 
\mbox{$ \psi^\ast \widetilde{I}^a (u\cdot f)= \widetilde{I}^a(\psi^\ast f).$}

In particular, for a hypersurface germ 
\mbox{$(F,z)\subset (\P^n\!,z)$} with isolated singularity
defined by $f$ we can introduce
\mbox{$\widetilde{I}^a(F,z):=\widetilde{I}^a(f)$} and
 \mbox{$\widetilde{X}^a(F,z):=V(\widetilde{I}^a(F,z))$}.
\end{lemma and definition}

\noindent {\it Proof.}
By the chain rule, we have \mbox{$j(g\circ\psi)=\psi^\ast(j(g))$} and,
obviously, \mbox{$j(u\!\!\:\cdot\!\!\: f)=j(f)$}.
\hfill $\Box$

\smallskip

\begin{lemma}
\label{generic}
Let \mbox{$f,g \in \K\{\underline{x}\}$} with $f$ an isolated 
singularity. Let \mbox{$\mathfrak{m}\subset \K\{\underline{x}\}$} be the
maximal ideal and let $j(f),j(g)$ denote the Tjurina ideals of $f,g$.
\begin{enumerate}
\itemsep2pt
\item[(a)] If \mbox{$j(g) \subset j(f)$} then \mbox{$f+tg$} is contact
equivalent to $f$ for almost all \mbox{$t\in \K$}.
\item[(b)] If \mbox{$j(g) \subset \mathfrak{m}\cdot j(f)$} 
then \mbox{$f+tg$} is contact
equivalent to $f$ for all \mbox{$t\in \K$}.
\item[(c)] If \mbox{$f+tg$} is contact equivalent to $f$ for sufficiently
small $t$, then
\mbox{$g \in \langle f\rangle + \mathfrak{m}\cdot j(f)$} 
\end{enumerate}
\end{lemma}

\noindent {\it Proof.} (a),(b) 
Set \mbox{$h_t:=f+tg$}. By assumption, there exists a
matrix \mbox{$A(\underline{x})=(a_{ij})_{i,j=0...n}$} such that
$$ \bigl( h_t, \tfrac{\partial h_t }{\partial x_1},\dots,
\tfrac{\partial h_t}{\partial x_n} \bigr) = 
\bigl( f, \tfrac{\partial f }{\partial x_1},\dots,
\tfrac{\partial f}{\partial x_n} \bigr) \cdot
\bigl(I+tA(\underline{x})\bigr)\,.$$ 
In Case (a) \mbox{det$\, 
\left(I\!+\!tA(0)\right)$} vanishes for at most \mbox{$n\!+\!1$} values of $t$,
while in Case (b) we have \mbox{det$\, 
\left(I\!+\!tA(0)\right)=1$} for all $t$ (since \mbox{$a_{ij}
\in \mathfrak{m}$}).
Since the Tjurina ideals $j(f)$ and $j(h_t)$ coincide if
\mbox{det$\, \left(I\!+\!tA(0)\right)\neq 0$}, (a) and (b) follow from
the  Theorem of Mather-Yau \cite{MatherYau}. (c) follows 
since $g$ is in the tangent space to the contact orbit, which is \mbox{$
\langle f\rangle + \mathfrak{m}\cdot j(f)$}. 
\hfill $\Box$

\smallskip
\begin{remark}
\label{Remark 1.9}
Since \mbox{$\mathfrak{m}_z^{\tau}\subset j(C,z)$} for 
\mbox{$\tau=\tau(C,z)$} the Tjurina number, Lemma \ref{generic}\:(b) says that
the local equation $f$ of $(C,z)$ 
is \mbox{$(\tau\!+\!1)$}-determined with respect to contact equivalence,
while  Lemma \ref{generic}\:(a) says that $f$ is 
\mbox{$\tau$}-deformation-determined. 
\end{remark}

\begin{lemma}
Let \mbox{$f\in \K\{\underline{x}\}$} be an isolated singularity. Then a
generic element \mbox{$g\in \widetilde{I}^a(f)$}
 is analytically equivalent to $f$ and
satisfies  \mbox{$\widetilde{I}^a(g)=\widetilde{I}^a(f)$}.

More precisely, let $d_0$ be the minimal degree of a polynomial defining
$\widetilde{I}^a(f)$. Then for any \mbox{$d\geq d_0$} the set of polynomials
in $\widetilde{I}^a(f)$ of degree \mbox{$\leq d$} which define 
$\widetilde{I}^a(f)$ is a Zariski-open dense subset. 
\end{lemma}

\noindent {\it Proof.} Let \mbox{$d\geq d_0$}. Then the polynomials 
\mbox{$g \in \widetilde{I}^a(f)$} of degree \mbox{$\leq d$} 
are parametrized by a finite dimensional vectorspace of
positive dimension. Since \mbox{$j(g)\subset j(f)$}, we have
\mbox{$\tau(g)\geq \tau(f)$} and equality holds exactly if \mbox{$
 j(g)=j(f)$}, that is, exactly if  \mbox{$\widetilde{I}^a(g)=
\widetilde{I}^a(f)$}. Now, the statement follows since the set of all $g$ with 
minimal possible Tjurina number \mbox{$\tau(g)=\tau(f)$} 
is a non-empty Zariski-open set.
\hfill $\Box$

\smallskip

\subsection{The Castelnuovo function of a zero-dimensional scheme in $\P^2$}
\label{sec:castelnuovo}

Let $X\subset \P^2$ be a zero-dimensional scheme. 

\begin{definition} The {\it Castelnuovo function} of $X$ is defined as 
$$
\kc_X:\;\Z_{\geq 0} \lra \Z_{\geq 0}\,,\quad
           d \longmapsto  h^1\bigl(\kj_{X/\P^2}(d-1)\bigr)-
           h^1\bigl(\kj_{X/\P^2}(d)\bigr)\,.
$$
\end{definition}

\noindent
In the following, we remind some basic properties of the Castelnuovo function,
which are obvious or can be proven by applying an elementary version of the 
so-called ``Horace method'' based on the exact sequence
$$ 0 \lra \kj_{X/\P^2}(d-1) \stackrel{\cdot L}{\lra} \kj_{X/\P^2}(d) \lra
\ko_L(d) \lra 0\,, 
$$ 
where $L$ denotes a generic line, 
respectively the corresponding exact cohomology sequence
$$ H^0\bigl(\kj_{X/\P^2}(d)\bigr) \lra
 H^0 \bigl(\ko_{L}(d)\bigr) \lra 
H^1\bigl(\kj_{X/\P^2}(d-1)\bigr)\lra H^1\bigl(\kj_{X/\P^2}(d)\bigr)\lra 0\,.$$
For the details, we refer to \cite{Davis}. 

\medskip\noindent
We introduce the notations
\begin{eqnarray*}
a(X)& = & \min\,\bigl\{\,d\in \Z\,\big|\,
h^0\bigl(\kj_{X/\P^2}(d)\bigr)>0\bigr\} \\
b(X)& = & \min\,\bigl\{\,d\in \Z\,\big|\, 
\big|H^0\bigl(\kj_{X/\P^2}(d)\bigr)\big|
\text{ has no fixed curve}\bigr\}\\ 
t(X)& = & \min\,\bigl\{\,d\in \Z\,\big|\, 
h^1\bigl(\kj_{X/\P^2}(d)\bigr)=0\bigr\}.
\end{eqnarray*}
Note that \mbox{$a(X)\leq b(X)\leq t(X)+1$}.
Let \mbox{$d\geq 0$} be an integer, then we have
\begin{enumerate} 
\itemsep2pt
\item[1.] \mbox{$\kc_Y(d) \leq \kc_X(d)$} for any subscheme \mbox{$Y\subset
X$}. 
\item[2.] \mbox{$\kc_X(0)+\ldots +\kc_X(d) =
  h^1\bigl(\kj_{X/\P^2}(-1)\bigr) -  h^1\bigl(\kj_{X/\P^2}(d)\bigr)= \deg X -
  h^1\bigl(\kj_{X/\P^2}(d)\bigr) $}. 
\item[3.] \mbox{$\kc_X(d)=0$} if and only if $d\geq t(X)+1$.
\item[4.] $\kc_X(d)\leq d+1$ with equality iff
  $h^0\bigl(\kj_{X/\P^2}(d)\bigr)=0$, that is, if $d\leq a(X)-1$.
\item[5.] if $d\geq a(X)$ then $\kc_X(d)\leq \kc_X(d-1)$.
\item[6.] if $b(X)\leq d \leq t(X)+1$ then $\kc_X(d) < \kc_X(d-1)$.

\begin{figure}[t]
  \centerline{\psfig{figure=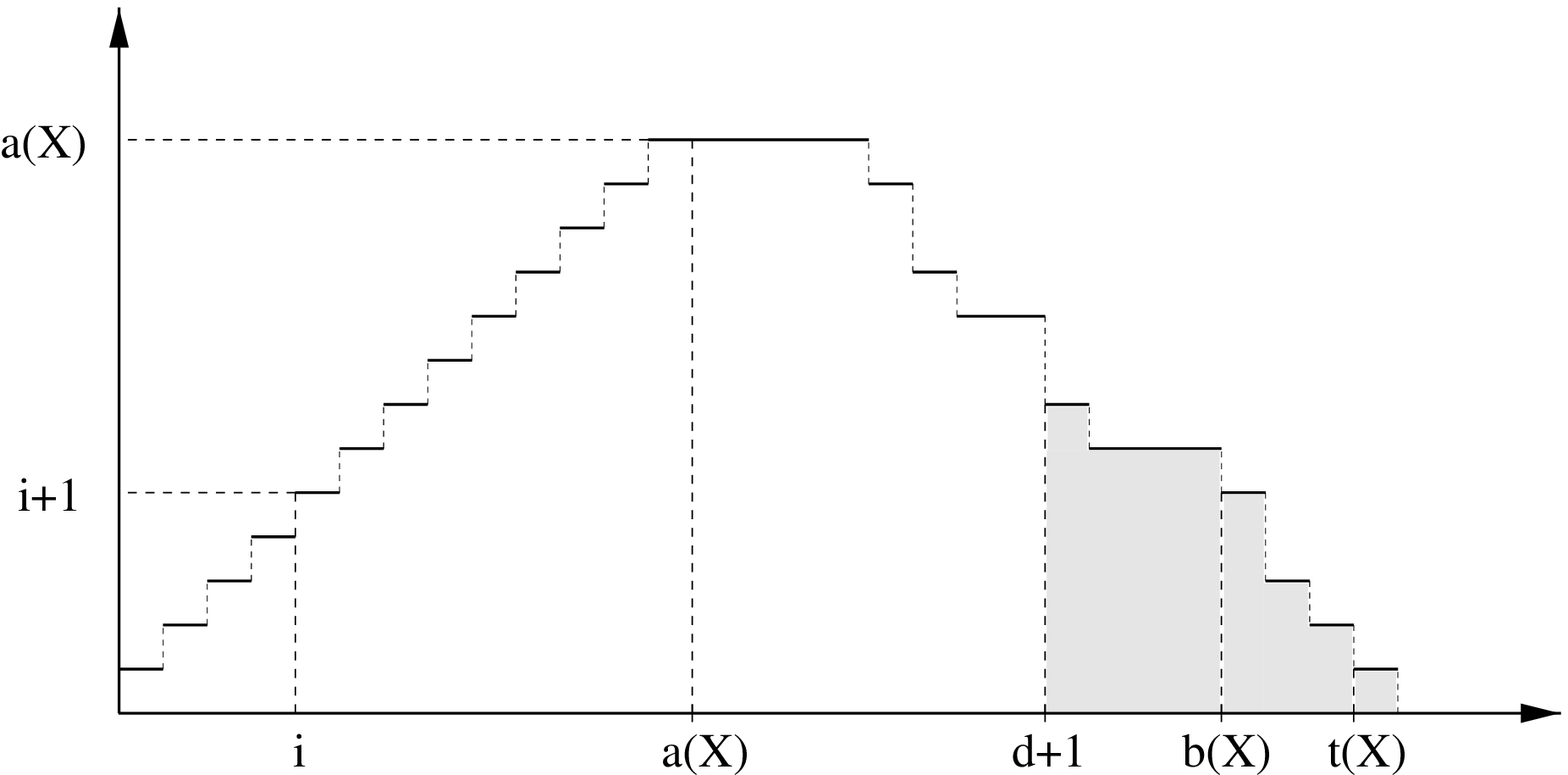,height=4.0cm,width=10.5cm}}
\caption{\label{fig 1} The graph of a Castelnuovo function (considered as a
  function on $\R_{\geq 0}$ given by \mbox{$\kc_X(t)=\kc_X([t])$}).
The content of the
  shaded region is $h^1\bigl(\kj_{X/\P^2}(d)\bigr)$.}
\end{figure}
\item[7.]
\begin{lemma}[Davis \cite{Davis}] 
\label{Davis}
Let \mbox{$X \subset \P^2$} be a zero-dimensional scheme and \mbox{$d_0\geq 
a(X)$} such that 
\mbox{$\kc_X(d_0)=\kc_X(d_0\!+\!1)$}. Then there exists a fixed curve $D$ of
degree $\kc_X(d_0)$ in the complete linear system
\mbox{$\big|H^0\bigl(\kj_{X/\P^2}(d_0)\bigr)\big|$} 
with the additional property that for each \mbox{$d\geq 0$} we have
\mbox{$ \kc_{X\cap D}(d)= \min \,\{\kc_X(d),\kc_X(d_0)\}.$}
\end{lemma}
\end{enumerate}

\begin{figure}
  \centerline{\psfig{figure=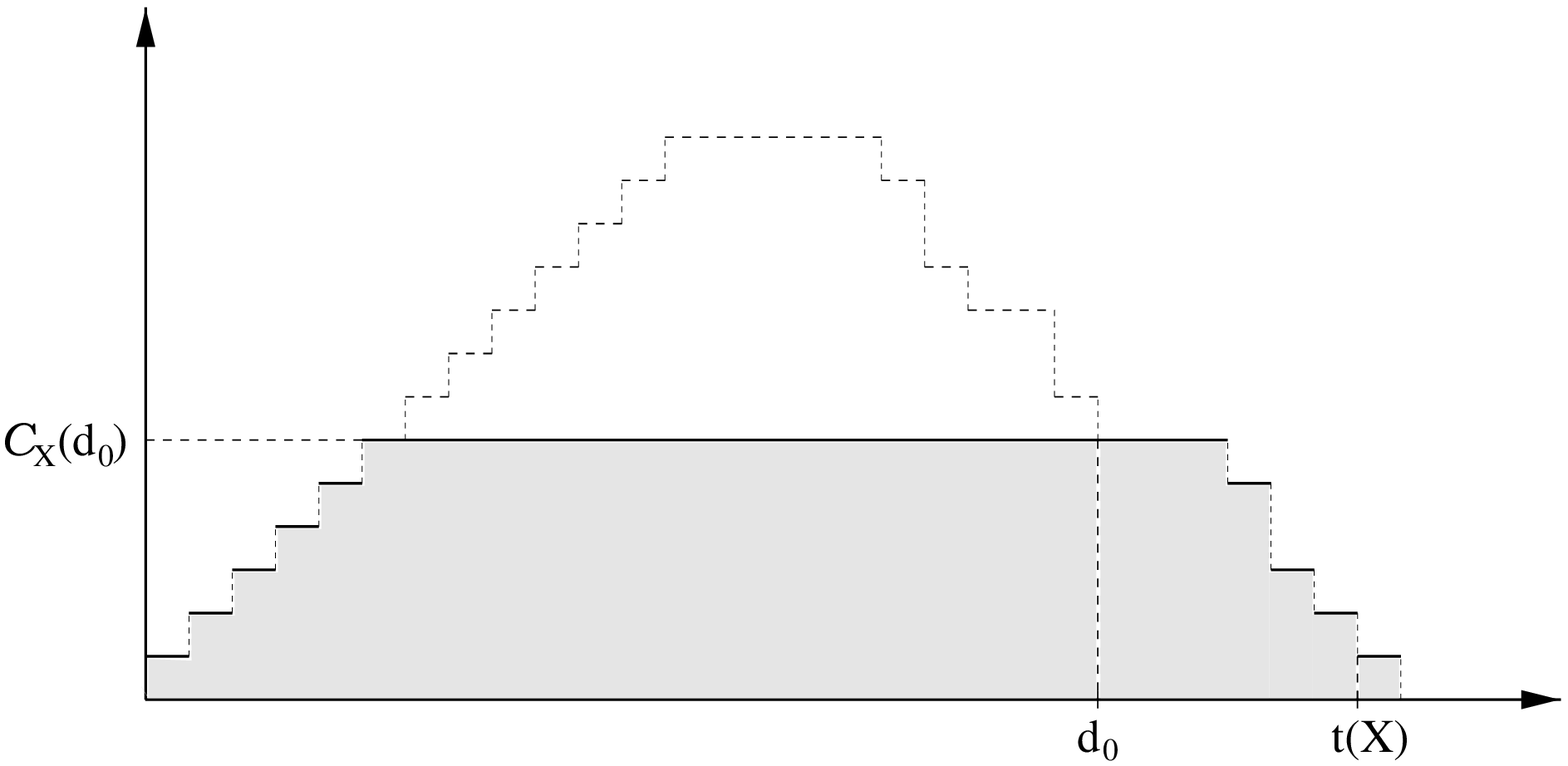,height=4.0cm,width=10.5cm}}
\caption{\label{fig 2} The graph of the Castelnuovo function $\kc_{X\cap
    D}$, where $D $ is the
  fixed curve in  $\big|H^0\bigl(\kj_{X/\P^2}(d_0)\bigr)\big|$ given by 
Lemma \ref{Davis}.
The content of the
  shaded region is \mbox{$\deg(X\cap D)$}.}
\end{figure}

\begin{definition}
We call a zero-dimensional scheme 
\mbox{$X\subset \P^2$} {\it decomposable} if there exists
a \mbox{$d_0>0$} such that
\mbox{$\kc_X(d_0\!-\!\!\:1)>\kc_X(d_0)=\kc_X(d_0\!+\!\!\:1)>0$}.
\end{definition}

\noindent
Finally, by B\'ezout's Theorem, we have
\begin{enumerate}
\item[8.] Let \mbox{$X=C_d\cap C_k$} be the intersection of two curves $C_d$
  and $C_k$ without common components. Moreover, let \mbox{$\deg C_d=d$},
  \mbox{$\deg C_k=k$} , \mbox{$k\leq d$}. Then 
  \mbox{$\kc_X(i)\leq k$} for each \mbox{$i\geq 0$} and
  \mbox{$\kc_X(d+k-i)=i-1$}  for any  \mbox{$i=1,\dots,k+1$}.
\end{enumerate}

\begin{figure}[t]
  \centerline{\psfig{figure=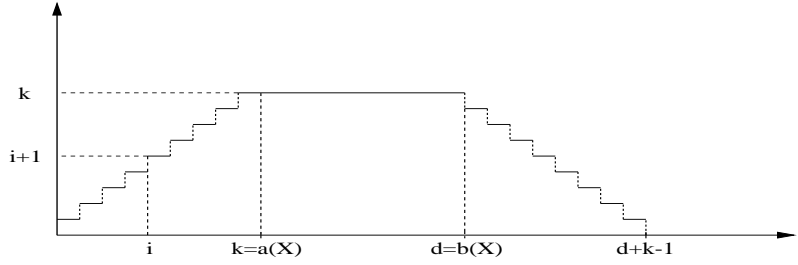,height=4.0cm,width=10.5cm}}
\caption{\label{fig 3} The graph of the Castelnuovo function for the complete
  intersection $X=C_d\cap C_k$.}
\end{figure}

\noindent
Considering these properties, it is not difficult to prove the following 
lemma, which is basically due to Barkats \cite{Barkats}.
\begin{lemma}
\label{Barkats} Let \mbox{$C_d\subset \P^2$} be an irreducible curve of degree
\mbox{$d>0$}, and \mbox{$X\subset C_d$}  a zero-dimensional scheme 
such that \mbox{$h^1\bigl(\kj_{X/\P^2}(d)\bigr) >0$}. Suppose moreover
\mbox{$ d> a(X)$}.
Then there exists a curve $C_k$ of degree \mbox{$k\geq 3$} such that the
scheme \mbox{$Y=C_k \cap X$} 
is non-de\-com\-po\-sable and satisfies
\begin{enumerate}
\itemsep0cm
\item \mbox{$h^1\bigl(\kj_{Y/\P^2}(d)\bigr)=h^1\bigl(\kj_{X/\P^2}(d)\bigr),$} 
\item \mbox{$k_0\cdot (d+3-k_0) \leq \deg Y,$} where \mbox{$k_0=\min \left\{
    k,\:\left[\tfrac{d+3}{2}\right]\right\}\,.$} 
\end{enumerate}
\end{lemma}

\noindent {\it Proof.}
{\it Case 1.} Suppose $X$ to be decomposable and 
let $d_0>0$ be maximal with the property \mbox{$\kc_X(d_0)=\kc_X(d_0\!+\!\!\:1)
>0.$} 

By Lemma \ref{Davis}, we obtain
the existence of a curve $C_k$ of degree \mbox{$k=\kc_X(d_0) <d$} such that
\mbox{$Y:=X\cap C_k$} is non-decomposable and 
\mbox{$ \kc_{Y}(i)= \min \{\kc_X(i),\kc_X(d_0)\}$} for each \mbox{$i\geq 0$}.
Remark that $Y$ is enclosed in the complete intersection \mbox{$C_d \cap
C_k$}, whence
\begin{equation}
\label{m-2}
1\:\leq\: \kc_Y(d+1)\:\leq \:\kc_{C_k \cap C_d} (d+1) \:=\: k-2
\,. 
\end{equation}
In particular, \mbox{$k\geq 3$} and 
$$ h^1\bigl(\kj_{Y/\P^2}(d)\bigr)\:= \textstyle{\sum\limits_{i=d+1}^{\infty}}
 \kc_Y(i) \:= 
\textstyle{\sum\limits_{i=d+1}^{\infty}}
 \kc_X(i)\:=\: h^1\bigl(\kj_{X/\P^2}(d)\bigr)\,.$$ 
Since $Y$ is non-decomposable
and \mbox{$\kc_Y(i)\leq k$} for each \mbox{$i\geq 0$}, we have
\begin{equation}
\label{k(d+2-k+1)}
\deg Y \: \geq \: \textstyle{\sum\limits_{i=0}^{d+1}} \kc_Y(i) 
\:\geq\: \textstyle{\sum\limits_{j=1}^{k_0}} 
\bigl(d+2-2(j\!-\!1)\bigr) \:=\: k_0(d+2-k_0+1) \,,
\end{equation}
whence the statement of the lemma.

\medskip \noindent
{\it Case 2.} If $X$ is a non-decomposable scheme, we can choose \mbox{$Y=X$}
and \mbox{$k:=a(X)<d$}. 
By the above reasoning we obtain again (\ref{m-2}) and (\ref{k(d+2-k+1)}). 
\hfill $\Box$

\smallskip
\begin{figure}[t]
  \centerline{\psfig{figure=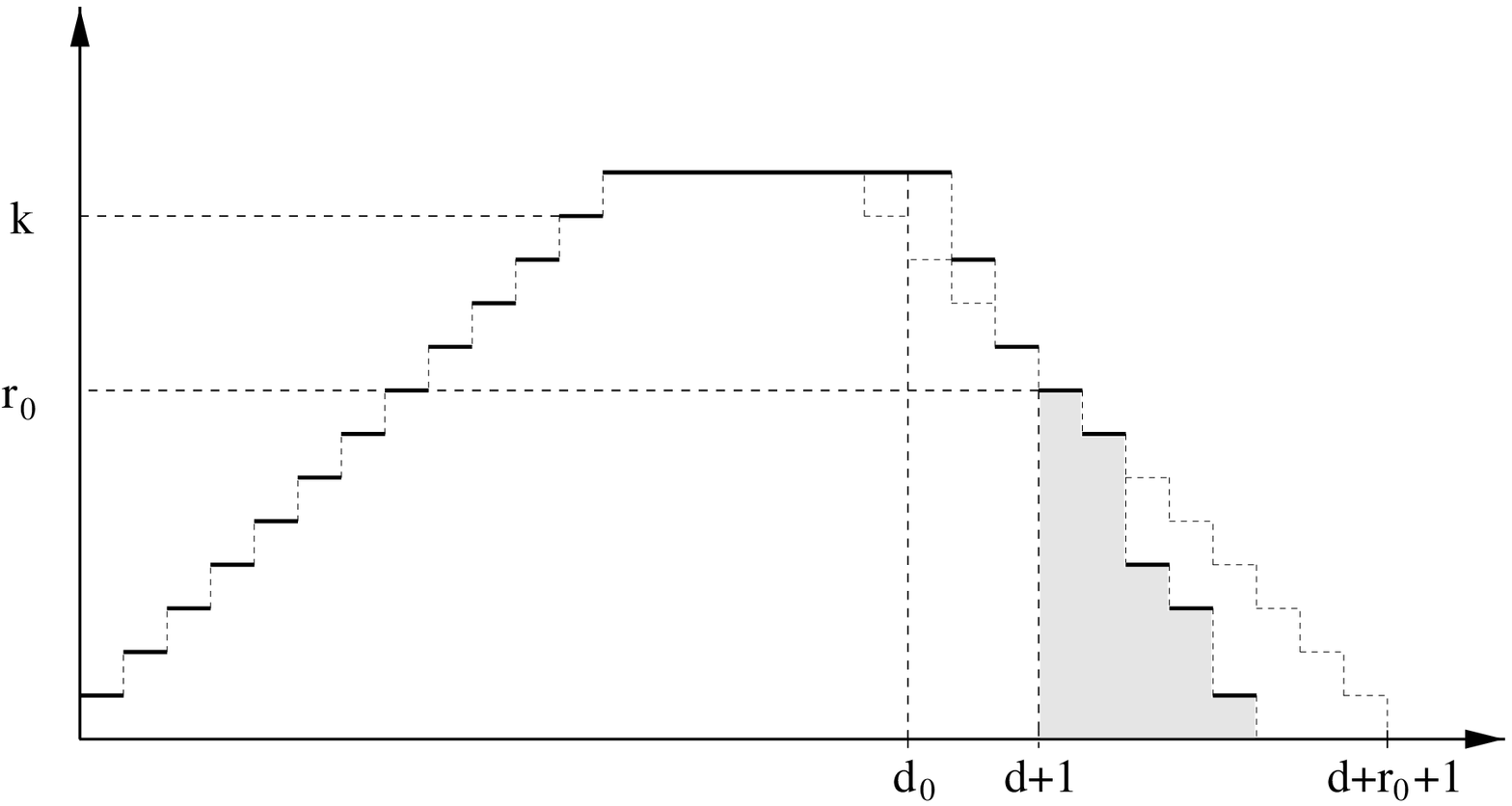,height=4.0cm,width=10.5cm}}
\caption{\label{fig 4} The graph of a Castelnuovo function $\kc_Y$. 
The content of the
  shaded region is $h^1\bigl( \kj_{Y/\P^2}(d)\bigr)$.}
\end{figure}

\begin{remark} 
\label{Barkats remark}
The zero-dimensional scheme $Y$ and the curve $C_k$ in Lemma
\ref{Barkats} satisfy 
\begin{equation}
\label{h1}
h^1\bigl(\kj_{Y/\P^2}(d)\bigr) \:=\,  \textstyle{\sum\limits_{i=d+1}^{\infty}}
 \kc_{Y} (i) 
 \:\leq \:\tfrac{r_0(r_0+1)}{2}\,,
\end{equation}
where \mbox{$r_0:=\kc_Y(d+1) \leq k-2$}
and (cf.\ Figure \ref{fig 4})
\begin{eqnarray}
\deg Y & = & h^1\bigl(\kj_{Y/\P^2}(d)\bigr) + \textstyle{\sum\limits_{i=0}^d}
 \kc_Y(i)\nonumber\\
&\geq&  h^1\bigl(\kj_{Y/\P^2}(d)\bigr) +
\textstyle{\sum\limits_{j=1}^{k_0}} 
(d+r_0+1-2(j\!-\!1))-\tfrac{r_0(r_0+1)}{2}  \label{degy} \\
& = & h^1\bigl(\kj_{Y/\P^2}(d)\bigr) + (d+2-k_0+r_0)\cdot k_0 -
\tfrac{r_0(r_0+1)}{2} \,.
\nonumber 
\end{eqnarray}
\end{remark}

\medskip
\section{Smoothness}
\setcounter{lemma}{0}
\setcounter{equation}{0}
\label{sec:smoothness}

\subsection{Equisingular and equianalytic families}
\label{sec:smoothness es}

Let \mbox{$S_1,\dots, S_r$} be topological (respectively analytic) types.
We recall that, by Proposition \ref{H0H1}, the variety \mbox{$\Vd$} is {\em
T-smooth} at \mbox{$C\in \Vd$} if and only if 
$$h^1\bigl(\kj_{\Xes(C)/\P^2}(d)\bigr)=0\,,\;\:\bigl(\text{respectively }\,
h^1\bigl(\kj_{\Xea(C)/\P^2}(d)\bigr)=0 \bigr)\,.$$ 
In order to formulate our results for the smoothness problem in a short way,
we first have to introduce new invariants for plane curve singularities. 
\begin{definition}
  Let \mbox{$(C,z)\subset (\P^2,z)$} be a reduced plane curve
singularity and 
$$\emptyset \neq X=X(C,z)\subset \Xea(C,z)$$
 be any zero-dimensional
scheme. Then we define for any
curve germ $(D,z)\subset (\P^2,z)$ without common
component with $(C,z)$
$$ \Delta (C,D;X):=\min\, \bigl\{\,(C,D)_z-
 \deg(D\cap X)\,, \; \deg(D\cap X) \,\bigr\}\,,$$
where $(C,D)_z$ denotes the
 local intersection number of the germs $(C,z)$ and $(D,z)$. Note that always
\mbox{$\Delta (C,D;X)\geq 1$} (cf.\
  Lemma \ref{Smoothness lemma} below). Hence, we can introduce 
\begin{equation}
\label{gamma es}
 \gamma\!\; \bigl(C;X\bigr) := \max_{(D,z)}
 \left\{ \tfrac{(\deg(D\cap X))^2}{\Delta (C,D;X)}
  + 2 \deg(D\cap X) + \Delta (C,D;X) \right\}\,,
\end{equation}
where the maximum is taken over all curve germs $(D,z)\subset (\P^2,z)$
that have no component in common with $(C,z)$. 
In particular, we introduce
$$ \gamma^{\text{es}} \!\:(C,z) := \gamma \!\;\bigl(C;\!\:\Xes(C,z)\bigr)\, 
\qquad \gamma^{\text{ea}}\!\: (C,z) := \gamma
\!\;\bigl(C;\!\:\Xea(C,z)\bigr)\,.$$ 
\end{definition}

\noindent
In Section \ref{sec:Section 5} we shall prove the following $h^1$-vanishing
theorem: 

\begin{proposition}
\label{Proposition 1}
Let \mbox{$C\subset \P^2$} 
be an irreducible curve of degree \mbox{$d\geq 6$} with $r$
singular points \mbox{$z_1,\dots,z_r$} and \mbox{$X_i\subset
  \Xea(C,z_i)$}, \mbox{$i=1,\dots,r$},
be any zero-dimensional schemes. Moreover, let $X$ be the (disjoint) union of
the schemes $X_1,\dots,X_r$. If
\begin{equation} 
\label{Smoothness condition 1} 
 \textstyle{\sum\limits_{i=1}^r} \gamma \!\;(C;X_i) \,<\, d^2+6d+8\,,
\end{equation}
then $h^1\bigl(\kj_{X/\P^2}(d)\bigr)$ vanishes.
\end{proposition}

\smallskip \noindent
As a corollary, we obtain our main smoothness result:

\begin{theorem}
\label{Smoothness Theorem}
Let \mbox{$C\subset \P^2$} be an irreducible curve of degree \mbox{$d\geq 6$} 
having $r$
singularities \mbox{$z_1,\dots,z_r$} of topological 
(respectively analytic) types $S_1,\dots,S_r$ as
its only singularities. Then 
\begin{enumerate}
\item[(a)] $\Vd$ is T-smooth at $C$ if
\begin{equation}
\label{Smoothness condition} 
\textstyle{\sum\limits_{i=1}^r} \gamma^{\text{es}}(C,z_i) < d^2+6d+8\quad
\Bigl(\text{respectively  }\, 
\textstyle{\sum\limits_{i=1}^r} \gamma^{\text{ea}}(C,z_i) < d^2+6d+8 \Bigr).
\end{equation}
\item[(b)] Under the condition
$$ \textstyle{\sum\limits_{i=1}^r} \gamma^{\text{ea}}(C,z_i) < d^2+6d+8 $$
the space of curves of degree $d$ is a joint versal deformation of
all singular points of the curve $C$. 
\end{enumerate}
\end{theorem}

\smallskip \noindent
In the following lemma we give general estimates for the invariants
$\gamma^{\text{es}} (C,z)$ (respectively $\gamma^{\text{ea}} (C,z)$) which 
show that Theorem \ref{Smoothness Theorem} improves the previously known
results (as stated above):

\begin{lemma}
\label{estim gamma}
For any reduced plane curve singularity \mbox{$(C,z)\subset (\P^2,z)$}, we can
estimate 
$$ \gamma^{\text{es}} (C,z) \:\leq (\tau^{\text{es}}(C,z)+1)^2\: \text{ and
  }\:  \gamma^{\text{ea}} (C,z) \:\leq (\tau (C,z)+1)^2\,,$$
where \mbox{$\tau(C,z)=\deg X^{\text{ea}} (C,z)$} denotes the Tjurina number,
  while \mbox{$\tau^{\text{es}}(C,z)=\deg X^{\text{es}} (C,z)$}
is the codimension of the $\mu$-const stratum in
  a versal deformation  base of $(C,z)$.  
\end{lemma}

\noindent {\it Proof.}
Let $(D,z)$ have no common component with
$(C,z)$, let \mbox{$X=X^{\text{es}}(C,z)$} (respectively 
\mbox{$X=X^{\text{ea}}(C,z)$}) and
\mbox{$\Delta= \Delta(C,D;X)$}.
There are two cases:

\smallskip
\noindent {\it Case 1.} \mbox{$\Delta=\deg D\cap X$}. Then, obviously,
$$  \tfrac{(\deg(D\cap X))^2}{\Delta} + 2 \deg(D\cap X) + \Delta \:=\:
 4 \cdot \deg (D\cap X) \:\leq \: 4 \cdot \deg X \:\leq \: (\deg X+1)^2.$$

\noindent {\it Case 2.} \mbox{$\Delta=(C,D)_z- \deg (D\cap X) < \deg
 (D\cap  X)$}, i.e.,  \mbox{$(C,D)_z<2\deg (D\cap X)$}. Then 
$$\tfrac{(\deg(D\cap X))^2}{\Delta} + 2 \deg(D\cap X) + \Delta \:=\: 
\tfrac{(C,D)_z^2}{(C,D)_z-\deg (D\cap X)}\,,$$
which is decreasing on 
\mbox{$\deg (D\cap X)+1\leq (C,D)_z \leq 2\deg (D\cap
X)-1$}. Consequently, it does not exceed
\mbox{$(\deg (D\cap X)+1)^2  \leq (\deg X+1)^2,$}
whence the statement.
\hfill $\Box$

\smallskip
\begin{examples}
\begin{enumerate}
\item[(a)] Let $(C,z)$ be an $A_k$-singularity (local equation
  \mbox{$x^2-y^{k+1}=0$}). Then we have 
 \mbox{$\gamma^{\text{es}} (C,z)=\gamma^{\text{ea}} (C,z)= (k+1)^2 =
(\tau^{\text{es}}(C,z)+1)^2$}.
\item[(b)] Let $(C,z)$ be a $D_4$-singularity (local equation
  \mbox{$x^3-y^3=0$}).  Then we obtain (cf.\ (\ref{gamma estimate}))  
\mbox{$\gamma^{\text{es}} (C,z)=\gamma^{\text{ea}} (C,z) = 18 < 25
=(\tau^{\text{es}}(C,z) +1)^2$}. 
\end{enumerate}
\end{examples}

\smallskip
\noindent
Applying the estimates from Lemma \ref{estim gamma}
to Theorem \ref{Smoothness Theorem}, we obtain in particular:

\begin{corollary} 
\label{cor a}
Let \mbox{$d\geq 6$} be an integer. Then
$\Vd$ is T-smooth at $C$ if 
$$
 \textstyle{\sum\limits_{i=1}^r}
 \left(\tau^{\text{es}}(C,z)+1\right)^2 < d^2+6d+8\quad
\Bigl(\text{respectively  }\, 
\textstyle{\sum\limits_{i=1}^r}
 \left(\tau (C,z)+1\right)^2 < d^2+6d+8 \Bigr)\,.
$$
\end{corollary}

\smallskip
\subsection{Families of curves with nodes and cusps}

Already for
families of curves with nodes and cusps, we obtain a slight improvement
against the previously known bounds (cf.\ \cite{GLS97,Shu97}).

\begin{corollary}
\label{Corollary 3.4}
The variety $\Vdi(n\cdot A_1,\,k\cdot A_2)$ of irreducible plane curves of
degree \mbox{$d\geq 6$} having $n$ nodes and $k$ cusps as its only
singularities is either empty or a smooth
variety of the expected dimension $d(d+3)/2-n-2k$ if 
$$ 4n+9k \: < \:d^2+6d+8\,.$$
\end{corollary}

\noindent {\it Proof.} This follows immediately from Theorem \ref{Smoothness
Theorem} (cf.\ also Corollary \ref{cor a}). 
\hfill $\Box$

\medskip
\subsection{Families of curves with ordinary singularities}

For families of curves with ordinary singularities (i.e., all local branches
are smooth and have different tangents) the new invariants pay off
drastically. We obtain a result which is not only 
asymptotically better than the previously known (cf.\ \cite{GLS97}), but even 
asymptotically proper.

\begin{corollary}
\label{Smoothn ord cor}
Let $\Vdi (m_1,\dots,m_r)$ be the variety of irreducible curves of degree
\mbox{$d\geq 6$} having $r$ ordinary multiple points of multiplicities
$m_1,\dots,m_r$, respectively, as only singularities. Then
\mbox{$\Vdi (m_1,\dots,m_r)$} is either empty or a smooth variety of the
expected 
dimension \mbox{$d(d+3)/2 -\sum_i \bigl(m_i(m_i\!\!\:+1)/2 -2\bigr)$} if
\begin{equation}
\label{Smoothness ord}
4\cdot \#(\text{nodes\/}) + 18 \cdot \# (\text{triple points\/}) +
\textstyle{\sum\limits_{m_i\geq 4}} 
\tfrac{16}{7} \cdot m_i^2 
\: < \: d^2+6d+8.  
\end{equation}
\end{corollary}

\noindent {\it Proof.} Let \mbox{$(C,z)\subset (\P^2,z)$}
be an ordinary $m$-fold point, then 
$\Ies:=\Ies(C,z)=j(C,z)+\mathfrak{m}_z^m$.
We shall show that 
\begin{equation}
\label{gamma estimate}
\max_{(D,z)}
 \left\{ \tfrac{\left(\deg(D\cap \Xes)+ 
\Delta^{\text{es}}(C,D)\right)^2}{\Delta^{\text{es}}(C,D)}\right\}
\:=\:
\gamma^{\text{es}} (C,z)\:\left\{
\begin{array}{ll}
=4 & \text{ if } m=2,\\ 
=18 & \text{ if } m=3,\\ 
\leq \tfrac{16}{7}m^2 & \text{ if } m \geq 4,
\end{array}
\right.
\end{equation}
whence (\ref{Smoothness ord}) implies (\ref{Smoothness condition}) and the
statement follows from Theorem \ref{Smoothness Theorem}.
Let $D=(D,z)$ be any plane curve germ of multiplicity $\mt D$ 
having no common component with $(C,z)$. 
As before, we have to consider two cases:

\smallskip \noindent
{\it Case 1.} \mbox{$\Delta^{\text{es}}(C,D)=\deg(D\cap \Xes)$}. Then 
$$  \gamma^{\text{es}} (C,z)\:= \:4 \deg(D\cap \Xes) \: \leq \: 4 \deg \Xes
\:=\: 2 m(m+1)-8 \,,$$
with equality if \mbox{$\mt D \geq m$}.

\smallskip\noindent
{\it Case 2.} \mbox{$\Delta^{\text{es}}(C,D)=(C,D)_z-\deg(D\cap
  \Xes)<\deg(D\cap \Xes)$}.
Note
that for fixed $\mt D$ and fixed $\deg(D\cap \Xes)$ the function
$$ \gamma\bigl((C,D)_z\bigr):= 
\tfrac{\left(\deg(D\cap \Xes)+ 
\Delta^{\text{es}}(C,D)\right)^2}{\Delta^{\text{es}}(C,D)}
\:=\:\tfrac{(C,D)_z^2}{ (C,D)_z-\deg(D\cap \Xes)}$$  
takes its maximum on \mbox{$\bigl[m\cdot \mt D,\, 2 \deg (D\cap \Xes)\bigr]$}
at \mbox{$(C,D)_z=m\cdot \mt D$}.
Hence, it is not difficult to see that it suffices to consider the cases 

\smallskip\noindent
{\it Case 2a.} \mbox{$m>\mt D=1$}, \mbox{$(C,D)_z=m$}. 
Then \mbox{$\deg (D\cap \Xes)=m-1$} and it follows that
\mbox{$\gamma\bigl((C,D)_z\bigr)= m^2$}. 

\smallskip \noindent
{\it Case 2b.} \mbox{$m>\mt D=2$}, \mbox{$(C,D)_z=2m$}.  
Then \mbox{$\deg (D\cap \Xes)=2(m-1)$}, which implies that
\mbox{$\gamma\bigl((C,D)_z\bigr)= 2m^2$}. 
 
\smallskip\noindent
{\it Case 2c.} \mbox{$m\geq\mt D\geq 3$}, \mbox{$m>3$}, 
\mbox{$(C,D)_z=m\cdot \mt D$}. 
Then
\begin{eqnarray*}
\deg (D\cap \Xes)&\leq & \dim_{\K} (\ko_{D,z} \big/ \mathfrak{m}_z^m) -1
\: = \: \tfrac{m(m+1)}{2}- 
\tfrac{ (m-\mt D)(m-\mt D+1)}{2} -1\\
&=&
m \cdot \mt D - \tfrac{(\mt D)^2 -\mt D+2}{2} \,,
\end{eqnarray*}
which implies that 
$$\gamma\bigl((C,D)_z\bigr)\:\leq\: \tfrac{2(m \cdot 
\mt D)^2}{(\mt D)^2 -\mt D+2} \: =\: 2m^2 \cdot
\tfrac{(\mt D)^2}{(\mt D)^2 -\mt D+2} \: \leq \: \tfrac{16}{7} \cdot
m^2\,.$$ 

\vspace{-20pt}
\hfill $\Box$

\medskip
\section{Irreducibility}
\setcounter{lemma}{0}
\setcounter{equation}{0}
\label{sec:irred}

\subsection{Equisingular and equianalytic families}
\label{sec:irred es}

Let \mbox{$S_1,\dots, S_r$} be topological (respectively analytic) types.
Moreover, let \mbox{$\nu'=\nu^s$} (resp.\ $\nu^a$) denote the 
deformation-determinacy as introduced in Section \ref{sec:topo} (respectively
\ref{sec:anal}) and \mbox{$\tau'=\tau^{\text{es}}$} (resp.\ 
\mbox{$\tau'=\tau$}) denote the codimension of the $\mu$-const stratum 
in the base of the semiuniversal deformation (respectively the Tjurina number).
Our main result on the irreducibility problem is:
\begin{theorem}
\label{Theorem 3}
If $d$ is a positive integer such that \mbox{$\max_i\,\nu'(S_i) \leq
  \tfrac{2}{5}d-1$}  and
\begin{eqnarray}
\textstyle{\sum\limits_{i=1}^r} \left(\nu'(S_i)+2\right)^2 
& < & \tfrac{9}{10} d^2
\,,  \label{(II)}\\
\tfrac{25}{2}
 \cdot \# (\text{nodes\/}) + 18 \cdot \# (\text{cusps\/})
+ \textstyle{\sum\limits_{\tau'(S_i)\geq 3}}
 \left(\tau'(S_i)+2\right)^2 & < & d^2 \label{(I)}
\end{eqnarray}
then $\Vd$ is irreducible or empty. 
\end{theorem}

\noindent
In particular, by Lemma \ref{Lemma 1.4} respectively Remark \ref{Remark 1.9},
we obtain the following, slightly weaker statement.
\begin{corollary}
If $d$ is a positive integer satisfying \mbox{$\,\max_i
  \tau'(S_i)\leq \tfrac{2}{5}d-1\,$} and
$$ \tfrac{25}{2}\cdot \# (\text{nodes\/})+18\cdot\#(
\text{cusps\/})+
\tfrac{10}{9} \cdot \textstyle{\sum\limits_{\tau'(S_i)\ge 3}}
 (\tau'(S_i)\!+\!2)^2 < d^2
$$
then $\Vd$ is irreducible or empty. 
\end{corollary}

\medskip \noindent
{\it Method of proof.} To be able to treat both, equisingular (es) and
equianalytic (ea),
families simultaneously, we introduce 
$$X(C):= \left\{
\begin{array}{cl}
X^s(C) & \text{ in the ``es''-case}\,,\\
X^a(C) & \text{ in the ``ea''-case}\,,
\end{array}
\right.
\;\;\,
\text{and}
\;\;\,
X'_{\text{fix}}(C):= \left\{
\begin{array}{cl}
\Xesf(C) & \text{ in the ``es''-case}\,,\\
\Xeaf(C) & \text{ in the ``ea''-case}\,.
\end{array}
\right. 
$$
Without restriction, we can assume that the types \mbox{$S_1,\dots,S_{r'}$},
\mbox{$r'\leq r$}, are pairwise distinct and that for any \mbox{$i=1,\dots,r'$}
the type $S_i$ occurs precisely $r_i$ times in \mbox{$S_1,\dots,S_r$}. We
introduce 
\begin{equation}\label{M}
\M=\M(S_1,\dots,S_r):=\textstyle{\prod\limits_{i=1}^{r'}
\text{Sym}^{r_i} (\P^2\!\times\! \kxf(S_i))}
\end{equation}
and consider the two morphisms
$$ 
\Phi_d:  \Vd \lra \text{Sym}^r\P^2\,, \quad C \longmapsto
(z_1\!+\!\ldots\!+\!z_r)\,,
$$
where 
\mbox{$(z_1\!+\!\ldots \!+\!z_r)$} is the unordered tuple of the singularities
of $C$ (cf.\ (\ref{Phi_d})), and
$$
 \Vd \underset{\text{dense}}{\supset} 
V \stackrel{\Psi_d}{\lra} \M\,,
\quad C \longmapsto 
\bigl[\bigl(z_i,\tau_{z_i0}(X(C,z_i))\bigr)\bigr]_{i=1,..,r}
$$
(cf.\ (\ref{Psi_d 2}), respectively (\ref{Psi_d 1})).
To obtain the irreducibility of \mbox{$\Vd$} or, equivalently,  of $V$, 
it suffices to prove that
the open subvariety
$$\Vr^{(2)}\,:=\, \left\{C \in V \,\big| \,
h^1\bigl(\kj_{X(C)/\P^2}(d)\bigr)=0\right\} \,\subset \, V$$
is dense and irreducible.

\medskip \noindent 
{\it Step 1.}  $\Vr^{(2)}$ is {\it irreducible}.

\smallskip \noindent
For any 
\mbox{$C\in \Vr^{(2)}$}, the fibre
\mbox{$\Psi_d^{-1} (\Psi_d(C))$} is the open dense subset $U$ of the linear
system 
\mbox{$\big|H^0\bigl(\kj_{X(C)/\P^2}(d)\bigr)\big|$}
 consisting of irreducible curves \mbox{$C'\in V$} with \mbox{$X(C')=X(C)$}.     
In particular, the fibres of $\Psi_d$ are smooth and equidimensional.
On the other hand, it follows from  Proposition
\ref{irreducible}, respectively Remark \ref{M irred},
that $\M$ is irreducible.
Hence, it suffices to show that $\Psi_d(
\Vr^{(2)})$ is dense in $\M$. This will be done in Section \ref{sec:Section 6}
(cf.\ Lemma \ref{dense}).

\medskip \noindent 
{\it Step 2.}  $\Vr^{(2)}$ is {\it dense} in \mbox{$V$}.

\smallskip \noindent
By Proposition \ref{H0H1}(e), we know that 
$$ \Vg \,:=\, \left\{ C\in  V \:\big| \: \text{Sing}\, C
  \text{ consists of points in general position}\right\}$$
is a dense subset of 
$$ \Vr^{(1)} \,:=\, \bigl\{ C \in V \,\big| \,
h^1\bigl(\kj_{X'_{\text{fix}}(C)/\P^2}(d)\bigr)=0\bigr\}\,.$$
Hence, it suffices to show that \mbox{$\Vg$} is a subset of $\Vr^{(2)}$ (this
will be done by applying a vanishing theorem for generic fat points, cf.\
Section \ref{sec:Section 6}) and that \mbox{$\Vr^{(1)}\subset V$} is dense. 
The latter statement takes the main part of Section \ref{sec:Section 6} and
will be proven by considering the Castelnuovo function 
associated to $X'_{\text{fix}}(C)$ (cf.\ Section \ref{sec:castelnuovo}). 
\hfill $\Box$

\smallskip
\medskip
\subsection{Families of curves with nodes and cusps}

Let \mbox{$\Vdi (n\!\!\:\cdot \!\!\:A_1,k\!\!\:\cdot \!\!\:A_2)$} 
be the variety of irreducible curves of degree $d$ having $n$ nodes and $k$
cusps as only singularities.  As an immediate corollary of Theorem 
\ref{Theorem 3}, we obtain:

\begin{corollary}
\label{nodes and cusps} Let \mbox{$d\geq 8$}. Then
the variety \mbox{$\Vdi (n\!\!\:\cdot \!\!\:A_1,k\!\!\:\cdot \!\!\:A_2)$}
 is irreducible or empty if
\begin{equation}
\label{n & c}
\tfrac{25}{2}\!\; n + 18\!\: k \: < \: d^2.
\end{equation}
\end{corollary}

\medskip
\subsection{Families with ordinary multiple points}

Let \mbox{$\Vdi (m_1,\dots,m_r)$} be the variety of irreducible curves of
degree $d$ having $r$ ordinary multiple points of multiplicities
$m_1,\dots,m_r$,  respectively, as only singularities. 
 \begin{corollary}
\label{ordinary singularities} Let \mbox{$\max \,m_i \leq \tfrac{2}{5}\!\;
d$}. Then \mbox{$\Vdi (m_1,\dots,m_r)$} is irreducible or empty if
\begin{equation}
\label{ord sin}
\tfrac{25}{2} \cdot \# (\text{nodes\/}) + 
\textstyle{\sum\limits_{m_i\geq 3}} \tfrac{m_i^2(m_i+1)^2}{4} \: < \: d^2. 
\end{equation}
\end{corollary}

\smallskip
\noindent {\it Proof.}  
This follows from Theorem \ref{Theorem 3}, since for an ordinary
$m_i$-tuple point $(C,z_i)$ we have
$$ \tau^{\text{es}}(C,z_i)+2=\deg  \Xesf(C,z_i) = \tfrac{m_i(m_i+1)}{2}\,, 
\qquad
\nu^s(C,z_i) = m_i-1\,.$$

\vspace{-18pt}
\hfill $\Box$

\smallskip
\subsection{Comments and Example}
We discuss here some aspects of the irreducibility problem
concerning the asymptotic properness of the results in Theorem 2 and
Corollary 3.1. To reach an asymp\-to\-ti\-cally proper sufficient
irreducibility condition one should try to improve the results obtained,
reducing singularity invariants in the left-hand side of the 
inequalities, or find examples
of reducible ESF with asymptotics of the singularity invariants
being as close as possible to that in sufficient conditions.

The classical problem of finding Zariski pairs, i.e., pairs of
curves of the same degree and with the same collection of singularities,
which have different fundamental groups of the complement in the plane,
has immediate relation to the problem discussed. Nori's theorem
\cite{No} states that {$\pi_1(\P^2\!\setminus\!\!\: C)=\Z/d\Z$}
for any curve \mbox{$C\in \Vd$} 
with
$$2\cdot\#(\mbox{nodes})+\!\textstyle{\sum\limits_{S_i\ne A_1}(\deg
  X^s(S_i) +\delta(S_i))\,<\,d^2}\ ,$$
where $X^{\text{es}}_{\text{fix}}$ is the zero-dimensional scheme defined in
Section \ref{sec:sec1.1} and $\delta(S_i)$ is the $\delta$-invariant. 
One can easily show that
the invariants in the left-hand side
 are $\le 3\mu$, hence any examples of Zariski
pairs must have asymptotics of singularity invariants as in the necessary 
condition for existence (\ref{??}), but not as in (\ref{(I)}).

%

The following proposition shows that
an equisingular family can have components of different dimensions,
whereas the fundamental groups of the complements of the curves are
the same.

\begin{proposition} \label{prop example}
Let $p,d$ be integers satisfying
\begin{equation}
p\ge 15,\quad 
6p\!\;<\!\;d\!\;\le \!\;12p-\tfrac{3}{2}-\sqrt{35p^2\!-\!\!\:15p\!\!\;+
\!\!\;\tfrac{1}{4}}\ .
\label{eqnew}\end{equation}
Then the family \mbox{$\Vdi(6p^2\!\cdot\!\!\: A_2)$} of irreducible 
curves of degree $d$ 
with $6p^2$ ordinary cusps has components of different dimensions. Moreover,
\mbox{$\pi_1(\P^2\!\!\:\setminus \!\!\;C)=\Z/d\Z$} for all \mbox{$\,C\in
  \Vdi(6p^2\!\cdot\!\!\: A_2)$}.
\end{proposition}

\begin{proof} Note that
(\ref{eqnew}) implies \mbox{$d^2>36p^2\!\!\:=6\!\!\:\cdot\!\!\: 6p^2$}. 
Hence, due to Nori's theorem (cf.\ \cite{No}), 
\mbox{$\pi_1(\P^2\!\!\:\setminus \!\!\; C)=\Z/d\Z$} for
all curves \mbox{$C\in \Vdi(6p^2\!\cdot\!\!\: A_2)$}.
We show that there are (at least) two different components of 
\mbox{$\Vdi(6p^2\!\cdot\!\!\: A_2)$}: by (\ref{eqnew}),
$$6p^2<\tfrac{(6p-1)(6p-2)+2}4<\tfrac{(d-1)(d-2)+2}4\ ,$$
and \cite{Sh5}, Theorem 3.3, gives the existence of a 
nonempty component $V'\!$ of
\mbox{$\Vdi(6p^2\!\cdot\!\!\: A_2)$} having the expected dimension
$$\dim V'=\tfrac{d(d+3)}2-12p^2$$
(the expected dimension in the construction of \cite{Sh5} 
follows from the $S$-transversality in \cite{Sh6}, Theorem 3.1). 

On the other hand, we construct a family $V''$ of bigger dimension: 
let $C_{2p}$, $C_{3p}$, $C'_{d-6p}$, $C''_{d-6p}$ 
be generic curves of degrees \mbox{$2p$},\,$3p,$
\mbox{$d\!\!\:-\!\!\:6p$}, \mbox{$d\!\!\:-\!\!\:6p$}, respectively. 
The curve 
\mbox{$C_d=C^3_{2p}C'_{d-6p}\!+C^2_{3p}C''_{d-6p}$} has degree $d$ and
$6p^2$ ordinary cusps as its only singularities, one at each intersection point
in \mbox{$C_{2p}\cap 
C_{3p}$}. Varying $C_{2p},C_{3p},C'_{d-6p},C''_{d-6p}$ in the spaces of curves
of degrees \mbox{$2p,\,3p,\,d\!\!\:-\!\!\:6p,\,d\!\!\:-\!\!\:6p$}, 
respectively, we obtain a subfamily $V''$ in
\mbox{$\Vdi(6p^2\!\cdot\!\!\: A_2)$}. Note that the equality
$$C_d\,=\,C^3_{2p}C'_{d-6p}\!\!\:+C^2_{3p}C''_{d-6p}\,=\,
\widehat{C}^3_{2p}
\widehat{C}'_{d-6p}\!\!\:+\widehat{C}^2_{3p}\widehat{C}''_{d-6p}\,=\,
\widehat{C}_d$$
with slightly deformed curves
$\widehat{C}_{2p},\widehat{C}_{3p},\widehat{C}'_{d-6p},\widehat{C}''_{d-6p}$ 
implies 
$$C_{2p}=\widehat{C}_{2p},\quad C_{3p}=\widehat{C}_{3p},\quad
C'_{d-6p}=\widehat{C}'_{d-6p},\quad C''_{d-6p}=\widehat{C}''_{d-6p}\, .$$
Indeed, if \mbox{$C_d=\widehat{C}_d$} 
then they have $6p^2$ common cuspidal points belonging to $C_{2p}$ 
and $\widehat{C}_{2p}$. Hence, by B{\'e}zout's
theorem, \mbox{$C_{2p}=\widehat{C}_{2p}$}. 
The tangent line to \mbox{$C_d=\widehat{C}_d$} at each cusp is
tangent to both, $C_{3p}$ and $\widehat{C}_{3p}$, that means, the intersection
number of \mbox{$C_{3p}$} and \mbox{$\widehat{C}_{3p}$} is at least
\mbox{$12p^2$},  whence 
\mbox{$C_{3p}=\widehat{C}_{3p}$}. We can conclude that
\mbox{$C_{2p}^3(C'_{d-6p}\!\!\:-\widehat{C}'_{d-6p})=
C_{3p}^2(\widehat{C}''_{d-6p}\!\!\:-C''_{d-6p})$} 
and, due to \mbox{$d\!\!\:-\!\!\:6p<2p$}, that 
\mbox{$C'_{d-6p}=\widehat{C}'_{d-6p}$},
\mbox{$C''_{d-6p}=\widehat{C}''_{d-6p}$}. Therefore, by (\ref{eqnew}),   
\begin{eqnarray*}
\dim V''&=&\textstyle{\frac{2p(2p+3)}2+\frac{3p(3p+3)}2+
2\cdot\frac{(d-6p)(d-6p+3)}2+1}\\
&=& \textstyle{\frac{d(d+3)}2-12p^2+\left(\frac{d^2}2-
d\left(12p-\frac{3}2\right)+ \frac{109p^2-21p+2}2\right)}
\:>\: \dim V'.
\end{eqnarray*}

\vspace{-17pt}
\end{proof}

\medskip
\smallskip

\section{Proof of Proposition \ref{Proposition 1}}
\setcounter{lemma}{0}
\setcounter{equation}{0}
\label{sec:Section 5}

\begin{lemma} 
\label{Smoothness lemma}
Let $(C,z)$ be a reduced plane curve singularity and let $I\subset
\mathfrak{m}_z \subset \ko_{\P^2,z}$ be an ideal containing the Tjurina ideal
$\Iea (C,z)$. Then for any $g\in I$ 
$$\dim_\K \ko_{\P^2,z}/ I \:<\:\dim_\K \ko_{\P^2,z}/\langle g,C\rangle \:=\:
(g,C)_z \,. $$
\end{lemma}

\noindent
{\it Proof.} cf.\ \cite{Shu97}, Lemma 4.1. \hfill $\Box$

\bigskip \noindent
Let $C$ be an irreducible curve of degree $d$ 
having precisely $r$ singularities $z_1,\dots,z_r$ and let 
$$X=X_1\cup \ldots\cup X_r\,, \quad
X_i\subset \Xea(C,z_i)\,,$$
\mbox{$i=1,\dots,r$}. 
Note that for any $i=1,\dots,r$ there exists a curve germ $(D,z_i)$ 
containing the scheme $X_i$ and satisfying \mbox{$\Delta(C,D;X_i)=\deg
  X_i$} (take any $(D,z_i)$ of sufficiently high multiplicity). 
Hence, we can estimate 
\begin{equation} 
\label{8888}
 \gamma(C;X_i) \: \geq \: \tfrac{\left(\deg X_i+\Delta(C,D;X_i)\right)^2}{
  \Delta(C,D;X_i)} \: = \: 4 \deg X_i\,. 
\end{equation}
In particular, by condition (\ref{Smoothness condition 1}) and since 
\mbox{$d\geq 6$}, we obtain
$$\deg X\,\leq\,\textstyle{\sum\limits_{i=1}^r}
 \tfrac{1}{4} \cdot \gamma (C;X_i) \,<\,
\tfrac{d^2+6d+8}{4} \,\leq \, \tfrac{d(d+1)}{2}\,,$$
whence \mbox{$d>a(X)=\min\,\big\{j\:\big| \: h^0\bigl(
\kj_{X/\P^2}(j)\bigr)>0\big\}$}. We want to show that
 \mbox{$h^1\bigl(\kj_{X/\P^2}(d)\bigr)$} vanishes.
Assume that this is not the case, that is, 
$$ h^1\bigl(\kj_{X/\P^2}(d)\bigr) >0\,.$$
Then  
Lemma \ref{Barkats} gives the existence of a curve $D$ of degree 
\mbox{$k\geq 3$}
such that \mbox{$Y=D\cap X$} is non-decomposable and satisfies 
\mbox{$h^1\bigl(\kj_{Y/\P^2}(d)\bigr)=h^1\bigl(\kj_{X/\P^2}(d)\bigr)
>0$}.  
Moreover, by (\ref{8888}) and (\ref{Smoothness condition 1}), we have 
\begin{equation}
\label{smoothness k<d/2}
\deg Y \:=\: \deg(X\cap D) \:< \: \tfrac{(d+3)^2}{4}-\tfrac{1}{4} \:\leq
\:\left[\tfrac{d+3}{2}\right]\cdot 
\left(d+3-\left[\tfrac{d+3}{2}\right]\right) \,.
\end{equation}
Hence, by Lemma \ref{Barkats}, \mbox{$k=k_0< \left[\tfrac{d+3}{2}\right]$} and
\begin{equation}
\label{Smoothness 1}
\deg Y \:\geq \: k \cdot (d+3-k)\,.
\end{equation}
Consequently, we can even estimate $k$ as
\begin{equation}
\label{Smoothness 2}
 k\:\leq\: \tfrac{d+3}{2} - \sqrt{\tfrac{(d+3)^2}{4} - \deg Y} \:=\:
\tfrac{2\cdot \deg Y}{d+3 + \sqrt{(d+3)^2 - 4\deg Y}}\,.
\end{equation}
On the other hand, let 
\mbox{$Y=Y_1 \cup \ldots \cup Y_s$}, \mbox{$\# Y:=s$},
be the decomposition of the zero-dimensional scheme $Y$ into its irreducible
components (without loss of generality, we may assume that $Y_i$ is supported
at $z_i$ for $i=1,\dots,s\leq r$). 
Note that, due to Lemma \ref{Smoothness lemma}, we have 
$$\deg Y_i \,\leq \,(C,D)_{z_i} - \Delta_i\,, \quad \Delta_i\geq 1\,, $$
which, together with  (\ref{Smoothness 1}) implies 
$$ k\cdot d \:\geq\: \textstyle{\sum\limits_{i=1}^s}  (C,D)_{z_i} \:\geq\: 
\deg Y +\textstyle{\sum\limits_{i=1}^s}
 \Delta_i \:\geq \: k \cdot (d+3-k)+\textstyle{\sum\limits_{i=1}^s}
 \Delta_i\,.$$
Thus, by (\ref{Smoothness 2}), we can estimate
$$\textstyle{\sum\limits_{i=1}^s}
 \Delta_i \:\leq \: k(k-3) \: < \: k^2 \:\leq \:  \Bigl(
\tfrac{2\cdot\deg Y }{d+3 + \sqrt{(d+3)^2 - 4\deg Y }}\Bigr)^2 \,.$$
In particular, applying the Cauchy inequality, we obtain
\begin{equation}
\label{Smoothness 3}
 \textstyle{\sum\limits_{i=1}^s}  \tfrac{(\deg Y_i)^2 }{\Delta_i} \: \geq \: 
\tfrac{(\deg
  Y)^2}{\Delta_1+\ldots+\Delta_s } \:>\: \tfrac{1}{4} \Bigl(1+ \sqrt{1 -
  \tfrac{4\deg Y}{(d+3)^2} } \:\Bigr)^2 \!\cdot (d+3)^2\,.
\end{equation}
Now, we introduce
$$ \alpha_Y:=\frac{\sum_{i=1}^s \frac{(\deg Y_i)^2}{\Delta_i}}{(d+3)^2}\,,\quad
\beta_Y:=\frac{\sum_{i=1}^s \frac{(\deg Y_i)^2}{\Delta_i}}{\deg Y} \,.$$
Then (\ref{Smoothness 3}) implies that
$$ \alpha_Y > \tfrac{1}{4} \cdot \Bigl(1 + \sqrt{1 -
  4\!\:\tfrac{\alpha_Y}{\beta_Y}}\:\Bigr)^2, \quad\text{ i.e., }\: \alpha_Y >
  \Bigl(\tfrac{\beta_Y}{\beta_Y+1}\Bigr)^2. $$
Finally, we have
\begin{eqnarray*}
(d+3)^2 & = & \tfrac{\beta_Y}{\alpha_Y} \cdot \deg Y \: < \: 
\bigl(1+\tfrac{1}{\beta_Y} \bigr)^2 \cdot \beta_Y \cdot \deg Y \:=\:
\bigl(\beta_Y + 2 + \tfrac{1}{\beta_Y}\bigr) \cdot \deg Y\\
& \leq & \textstyle{\sum\limits_{i=1}^s} \Bigl( \tfrac{(\deg Y_i)^2}{\Delta_i}
+ 2\cdot 
  \deg Y_i + \Delta_i \Bigr) \: \leq \: \textstyle{\sum\limits_{i=1}^s}
 \gamma (C;X_i)\,,
\end{eqnarray*}
which contradicts (\ref{Smoothness condition 1}).
\hfill $\Box$

\medskip
\section{Proof of Theorem \ref{Theorem 3}}
\setcounter{lemma}{0}
\setcounter{equation}{0}
\label{sec:Section 6}

\noindent
In this section, we complete the proof of Theorem \ref{Theorem 3}. To do so,
using the notations introduced in Section
\ref{sec:irred es}, we shall prove the following lemmas:

\begin{lemma} 
\label{dense}
If $\Vr^{(2)}$ is non-empty 
then \mbox{$\Psi_d (\Vr^{(2)})$} is dense in $\M$.
\end{lemma}

\begin{lemma} 
\label{Lemma 6.4}
Let \mbox{$C\in \Vd$} be a curve that has its singularities 
in generic position \mbox{$z_1,\dots,z_r$}. If \mbox{$ \,2d> 5 \cdot\max_i
\nu'(C,z_i) + 4\,$} and
\begin{equation} 
\label{Xu cond}
\tfrac{9}{10}\cdot
 (d+3)^2 > \textstyle{\sum\limits_{i=1}^r} \bigl(\nu'(C,z_i)+2\bigr)^2,
\end{equation}
then \mbox{$h^1\bigl(\kj_{X(C)/\P^2}(d)\bigr)$} vanishes, that is,
$\Vg$ is a subset of $\Vr^{(2)}$.
\end{lemma}

\begin{lemma} 
\label{Lemma 6.3}
Let \mbox{$d\geq 8$} be an integer and \mbox{$C\in \Vd$} such
  that 
\begin{eqnarray}
\label{cond0}
d^2+6d+8 & > & 4 \,\deg X'_{\text{fix}}(C)\,,\\
\label{cond1} 
d^2 & > & \textstyle{\sum\limits_{i=1}^r} (\deg X'_{\text{fix}}(C,z_i))^2,\\
2\cdot(d+3)^2 & > & \textstyle{\sum\limits_{i=1}^r}
 (\deg X'_{\text{fix}}(C,z_i)+2)^2,
\label{cond2}\\ 
\tfrac{9}{10} \cdot d^2 & > & \textstyle{\sum\limits_{i=1}^r} 
\max \, \bigr\{\, \bigl( \deg D \cap X'_{\text{fix}}(C,z_i) \bigr)^2 \,
\big| \: D \text{ a smooth curve}\,\bigr\}, \label{cond3}\\
(d-1)^2 & > & \textstyle{\sum\limits_{i=1}^r} 
\max\,\left(
\begin{array}{c}
\bigl\{\,\bigl(\deg D\cap X'_{\text{fix}}(C,z_i) \bigr)^2\,\big|\:
D \text{ a smooth curve}\,\bigr\}\\
\cup \:\bigl\{\,\tfrac{1}{2} \cdot
\bigl(\deg X'_{\text{fix}}(C,z_i)\bigr)^2\,\bigr\}
\end{array}
\right), \label{cond4}\\ 
\tfrac{16}{15} \cdot (d+3)^2 & > & \textstyle{\sum\limits_{i=1}^r }
\max\,\left(
 \begin{array}{c}
\bigl\{\,\bigl(\deg D\cap
X'_{\text{fix}}(C,z_i) +\tfrac{16}{15}\bigr)^2\,\big|\: 
D \text{ a smooth curve} \,\bigr\}\\
\cup \:\bigl\{\,\tfrac{1}{2} \cdot 
\bigl(\deg
  X'_{\text{fix}}(C,z_i)+\textstyle{\frac{32}{15}}\bigr)^2\,\bigr\}
\end{array}
\right).  \label{cond5} 
\end{eqnarray}
Then $\Vr^{(1)}$ is dense in $V$ ,i.e.,
\mbox{$h^1 \bigl( \kj_{X'_{\text{fix}}(C)/ \P^2}(d)\bigr) =0$} for
generic \mbox{$C\in V$}.
\end{lemma}

\begin{remark}
Note that for any reduced plane curve singularity \mbox{$(C,z)\subset
  (\P^2,z)$} and any smooth curve germ $D$ at $z$  we have 
$$ \deg X'_{\text{fix}}(C,z)\:=\:\tau'(C,z)+2 \:\geq\:
\nu'(C,z)+2\,, \quad \deg\bigl(D\cap X'_{\text{fix}}(C,z)\bigr) \: \leq \:
\nu'(C,z)+1\,. $$
For instance, in the case of nodes and cusps, we have
$$ \deg  X'_{\text{fix}}(C,z) = \left\{
\begin{array}{ll}
 3 & \text{ for a node},\\
 4 & \text{ for a cusp},
\end{array}
\right.
\qquad
\nu'(C,z) =
\left\{
\begin{array}{ll}
 1 & \text{ for a node},\\
 2 & \text{ for a cusp},
\end{array}
\right.
$$
$$
\max \, \bigl\{\!\:\deg  (D \cap X'_{\text{fix}}(C,z)) \:\big|\: D \text{
  smooth}\,\bigr\} \, =\,
\left\{
\begin{array}{ll}
 2 & \text{ for a node},\\
 3 & \text{ for a cusp}.
\end{array}
\right.
$$
Hence, it is not difficult to see that the conditions (\ref{(I)}) and
(\ref{(II)}) imply (\ref{cond0})--(\ref{cond5}).
\end{remark}

\medskip
\noindent {\it Proof of Lemma \ref{dense}.} By Sections \ref{sec:topo},
\ref{sec:anal}, 
we know that for any \mbox{$i=1,\dots,r$} there exists an 
\mbox{$m_i$} such that the schemes \mbox{$X(C,z_i)$}, 
depend only on the \mbox{$(m_i\!-\!1)$-jet} of the equation of
$(C,z_i)$. Hence, 
for \mbox{$d_0 \geq \max\, m_i$} the morphism $\Psi_{d_0}$ is dominant. 
Moreover, we can assume $d_0$ to be sufficiently large such that 
$h^1\bigl(\kj_{X(C)/\P^2}(d_0)\bigr)$ vanishes. Hence,
\begin{eqnarray*}
\dim \,\M &=&  \dim \Psi_{d_0} (V)\:=\:
\dim V_{d_0}(S_1,\dots,S_r) - h^0\bigl(\kj_{X(C)/\P^2}(d_0)\bigr)+1\\
& = & \dim V_{d_0}(S_1,\dots,S_r) - \tfrac{d_0(d_0+3)}{2} +\deg X(C) \,.
\end{eqnarray*}
On the other hand, 
let \mbox{$C\in \Vr^{(2)}$}. Then    
the vanishing of $h^1\bigl(\kj_{X(C)/\P^2}(d)\bigr)$ implies in particular 
the T-smoothness of 
$V_d(S_1,\dots,S_r)$ at $C$ (cf.\ Proposition \ref{H0H1}\:(c)). Hence, 
as an open subvariety, $\Vr^{(2)}$ is also smooth at $C$ 
of the expected codimension
$$\tfrac{d(d\!\!\;+\!\!\;3)}{2} - \dim \Vr^{(2)}\:=\:
\tfrac{d_0(d_0\!\!\:+\!\!\;3)}{2} - \dim V_{d_0}(S_1,\dots,S_r)\,,$$
that is,
\begin{eqnarray*}
\dim \Psi_{d} (\Vr^{(2)}) & = & 
\dim \Vr^{(2)} - h^0\bigl(\kj_{X(C)/\P^2}(d)\bigr) +1\: = \:
 \dim \Vr^{(2)} -\tfrac{d(d+3)}{2} + \deg X(C) \\
& = & \dim V_{d_0}(S_1,\dots,S_r) -\tfrac{d_0(d_0+3)}{2} + \deg X(C)
\: =\: \dim \M\,,
\end{eqnarray*}
whence the statement.
\hfill $\Box$

\bigskip \noindent
{\it Proof of Lemma \ref{Lemma 6.4}.}
Let \mbox{$i\in \{1,\dots,r\}$} and \mbox{$\nu_i=\nu'(C,z_i)$}.
By definition of $\nu_i$, 
the scheme $X(C,z_i)$ is contained in the ordinary fat point given by
the ideal $\mathfrak{m}_{z_i}^{\nu_i+1}$. Hence it suffices to show that
\mbox{$h^1\bigl(\kj_{Y(\nu_1+1,\dots,\nu_r+1)/\P^2}(d)\bigr)=0\,,$} 
where \mbox{$Y(\nu_1\!+\!1,\dots,\nu_r\!+\!1)$}
 is the zero-dimensional scheme of $r$ ordinary fat points
of multiplicities \mbox{$\nu_1\!+\!1,\dots, \nu_r\!+\!1$} in general position.
Now, the statement follows from \cite{GengXu95}, Theorem 3.
\hfill $\Box$

\bigskip \noindent
{\it Proof of Lemma \ref{Lemma 6.3}.}
Assume $V$ has an irreducible component
\mbox{$V^\ast \!\subset V\backslash \Vr^{(1)}$}, that is, the 
generic element $C$ of $V^\ast$ satisfies 
$$h^1 \bigl(\kj_{X'_{\text{fix}}(C)/\P^2}(d) \bigr) >0\,.$$ 
We denote by \mbox{$\Sigma^\ast\subset \text{Sym}^r \P^2=:\Sigma$} the closure
of $\Phi_d(V^\ast)$. 

\setlength{\unitlength}{1cm}
\begin{picture}(13,3.4)
\put(5.1,0.25){$\Sigma^\ast$}
\put(6.5,0.25){$\subset $}
\put(6.3,0){$\scriptstyle{\text{closed}}$}
\put(7.6,0.25){$\text{Sym}^r \P^2=:\Sigma$}
\put(4.23,0.8){$\scriptstyle{\text{dense}}$}
\put(5.1,0.8){$\cap$}
\put(4.7,1.3){$\Phi_d(V^\ast)$}
\put(5.2,2.6){\vector(0,-1){0.8}}
\put(4.6,2.2){$\textstyle{\Phi_d}$}
\put(8.2,2.6){\vector(0,-1){1.9}}
\put(8.4,1.7){$\textstyle{\Phi_d}$}
\put(8.1,2.8){$V$}
\put(7.5,2.8){$\subset $}
\put(6.2,2.8){$V\backslash \Vr^{(1)}$}
\put(5.65,2.8){$\subset $}
\put(5.0,2.8){$V^\ast$}
\end{picture}

\noindent
Recall that the dimension of 
\mbox{$\Phi_d^{-1}\bigl(\Phi_d(C)\bigr)$} at $C$ is just the dimension of
\mbox{$V_{d,\text{fix}} (S_1,\dots,S_r)$} at $C$, that is, by Proposition
\ref{H0H1}\:(b),
$$ \dim \Phi_d^{-1}\bigl(\Phi_d(C)\bigr) \:\leq \:
h^0\bigl(\kj_{X'_{\text{fix}}(C)/\P^2}(d)\bigr) -1\,.$$
To obtain the statement of  Lemma \ref{Lemma 6.3}, it
suffices to show that under the given (numerical) conditions we would have
\begin{equation}
\label{codim}
h^1\bigl(\kj_{X'_{\text{fix}}(C)/\P^2}(d)\bigr) < \:\text{codim}_\Sigma
\Sigma^\ast , 
\end{equation}
because this would imply that 
\begin{eqnarray*}
\dim V^\ast &\leq & \dim \Sigma^\ast +
h^0\bigl(\kj_{X'_{\text{fix}}(C)/\P^2}(d)\bigr) -1 \\
&<& \dim \Sigma + h^0\bigl(\kj_{X'_{\text{fix}}(C)/\P^2}(d)\bigr) - 
h^1\bigl(\kj_{X'_{\text{fix}}(C)/\P^2}(d)\bigr) -1 \,=\, \dim \Vr^{(1)} \,,
\end{eqnarray*}
whence a contradiction (any component of $V$ has at least the expected 
dimension $\dim \Vr^{(1)}$).

\medskip \noindent
{\it Step 1.} For \mbox{$d\geq 6$} the condition (\ref{cond0}) implies in
particular that  
\mbox{$\deg X'_{\text{fix}}(C) \leq  d(d\!\!\:+\!\!\:1)/2$},
whence  
\mbox{$ d > a(X'_{\text{fix}}(C)) = \min\,\bigl\{i\,\big|\,
h^0\bigl(\kj_{X'_{\text{fix}}(C)/\P^2}(i)\bigr)>0\bigr\}$}. 
By Lemma \ref{Barkats}, we obtain the existence of a curve $C_k$ of
degree \mbox{$k\geq 3$} such that the subscheme 
\mbox{$Y=C_k\cap X'_{\text{fix}}(C) \subset C_k\cap C$}
is non-decomposable with
\begin{equation}
\label{Bedingung 1a}
h^1\bigl(\kj_{Y/\P^2}(d)\bigr)\, =\,
h^1\bigl(\kj_{X'_{\text{fix}}(C)/\P^2}(d)\bigr) \:\leq  
\:\tfrac{r_0(r_0+1)}{2}\,, 
\end{equation} 
where \mbox{$1 \leq r_0:=\kc_{X'_{\text{fix}}(C)}(d+1) \leq k-2$} (cf.\ Remark 
\ref{Barkats remark}). 
Since, by (\ref{cond0}), we suppose additionally that 
\begin{equation}
\label{k<d/2}
\deg Y \:< \: \left[\tfrac{d+3}{2}\right]\cdot
\left(d+3-\left[\tfrac{d+3}{2}\right]\right) ,
\end{equation}
we have \mbox{$k< \left[\tfrac{d+3}{2}\right]$} 
and (cf.\ Lemma \ref{Barkats} and Remark \ref{Barkats remark}) 
\begin{equation}
\label{Bedingung 1}
\deg Y \:\geq \: \max\left\{ k \cdot (d\!+\!3\!-\!k)\,,\;\:
k\cdot(d\!+\!2\!+\!r_0\!-\!k)+ 
h^1\bigl(\kj_{Y/\P^2}(d)\bigr)- \tfrac{r_0(r_0+1)}{2} \right\}\,.
\end{equation}
Now, we can estimate the codimension of $\Sigma^\ast$ in $\Sigma$.
Given the curve
$C_k$, the number of conditions on $X'_{\text{fix}}(C)$ imposed by fixing the support of the
subscheme \mbox{$Y=C_k\cap X'_{\text{fix}}(C)$} on $C_k$ respectively its singular locus
is at least
\mbox{$\# Y$} if $C_k$ is non-reduced and at least 
\mbox{$\# Y + \# (Y|_{\text{Sing}\,C_k})$} if $C_k$ is a reduced curve.
On the other hand, the dimension of the variety of reduced (respectively non
reduced) curves $C_k$ of degree $k$ is given by
\mbox{$h^0\bigl(\ko_{\P^2}(k)\bigr)-1$} (respectively 
\mbox{$h^0\bigl(\ko_{\P^2}(k\!\!\:-\!\!\:2)\bigr)+2$}).
Thus, in place of (\ref{codim}), it suffices to show that
\begin{equation}
\label{codim2}
 h^1\bigl(\kj_{Y/\P^2}(d)\bigr)\: <\: \min\, \left\{
\# Y - \tfrac{k^2-k}{2}-2\,,\; \# Y + \# (Y|_{\text{Sing}\,C_k})- 
\tfrac{k(k+3)}{2}\right\} \,.
\end{equation}

\medskip
\noindent
{\it Step 2.}
Recall that we have \mbox{$k\geq 3$} and, by (\ref{Bedingung 1a}), 
\begin{equation}
\label{h1neu} 
h^1\bigl(\kj_{X'_{\text{fix}}(C)/\P^2}(d)\bigr) \:=\:
h^1\bigl(\kj_{Y/\P^2}(d)\bigr) \:\leq\: 
\tfrac{(k-2)(k-1)}{2}\,. 
\end{equation}

\smallskip\noindent
{\it Step 2a.} Assume \mbox{$\,h:=h^1\bigl(\kj_{Y/\P^2}(d)\bigr)=
\tfrac{(k-2)(k-1)}{2}$}.

\smallskip\noindent
Note that this implies that the Castelnuovo functions of $Y$ and 
\mbox{$C_k\cap C$} coincide, in particular we have \mbox{$\deg Y = kd$},   
i.e., \mbox{$Y=C_k \cap C$}.
In this case the condition (\ref{codim2}) is satisfied whenever
\begin{equation}
\label{Step 2a}
 0 \: < \:\min \,\left\{\# Y - k^2+2k-3\,,\; \# Y + 
\# (Y|_{\text{Sing}\,C_k})- 
k^2-1\right\} \,.
\end{equation}
Now, we have to consider two cases

\smallskip\noindent
{\it Case 1:} \mbox{$\# (Y|_{\text{Sing}\,C_k})\geq 1$}. 
Then the right-hand side
is bounded from below by $\# Y-k^2 = \# Y - (\deg Y)^2/d^2,$ 
whence, due to the Cauchy inequality, it suffices to have 
\begin{equation}
\label{cond1neu} 
d^2 \: > \: \textstyle{\sum\limits_{i=1}^r}
 \deg (\underbrace{X'_{\text{fix}}(C,z_i)\cap C_k}_{
\displaystyle{=:Y_i}})^2, 
\end{equation}
which is implied by (\ref{cond1}).

\smallskip\noindent
{\it Case 2:} \mbox{$\# (Y|_{\text{Sing}\,C_k})= 0$}. 
Then, as \mbox{$k\geq 3$}, the right-hand side is bounded from below by 
$\# Y-k^2-1\geq \# Y - \tfrac{10}{9} k^2$,
whence (\ref{Step 2a}) holds whenever
$$ \textstyle{\sum\limits_{i=1}^r}
 (\deg Y_i)^2 \: < \: \tfrac{9}{10}\, d^2 \;\text{ with }\;
\deg Y_i \: \leq \: \max\,\left\{\deg D \cap X'_{\text{fix}}(C,z_i) \mid D 
\text{ smooth}\right\}\,, 
$$
which is a consequence of (\ref{cond3}).

\smallskip
\medskip\noindent
{\it Step 2b.} Assume 
\mbox{$h=h^1\bigl(\kj_{Y/\P^2}(d)\bigr)<\tfrac{(k-2)(k-1)}{2}\,,$}  
in particular \mbox{$k\geq 4$}.

\smallskip \noindent
As we have seen in (\ref{codim2}), it suffices to show that
\begin{equation}
\label{Step2b}
\max_{k,h}\, \Bigl\{ h+\underbrace{\tfrac{k^2-k}{2}+2}_{\displaystyle
=:p_1(k)}\Bigr\}\: 
< \:\# Y \;\, \text{ and }\;\, 
\max_{k,h}\, \Bigl\{ h +\!\!\underbrace{\tfrac{k(k+3)}{2}}_{\displaystyle
=:p_2(k)}\Bigr\}\: 
< \:\# Y + \# (Y|_{\text{Sing}\,C_k}) \,.
\end{equation}
We introduce
$$\rho_j\::=\: \min \,\bigg\{\tfrac{(\deg Y)^2}{p_j(k)+h} \;\Big|\;
\begin{array}{l}
1 \leq h \leq \min\,\bigl\{ \tfrac{r_0(r_0+1)}{2},\,\tfrac{k(k-3)}{2}\bigr\}\\
4 \leq k\,, \; 1\leq r_0 \leq k-2 
\end{array}
\bigg\} \,, \qquad j=1,2\,.
$$
By the Cauchy inequality, it follows that (\ref{Step2b}) holds whenever
\begin{equation}
\label{beta}
\textstyle{\sum\limits_{i=1}^s} (\deg Y_i)^2
\:<\:\rho_1 \;\, \text{ and }  \textstyle{\sum\limits_{z_i \not\in
    \text{Sing}\,C_k} } 
(\deg Y_i)^2 \,+ \textstyle{\sum\limits_{z_i \in \text{Sing}\,C_k}} 
\tfrac{(\deg Y_i)^2}{2}\:<\:\rho_2 \,.
\end{equation}
It remains to estimate $\rho_1$ and $\rho_2$ as functions in $d$. 
By (\ref{Bedingung 1}), 
we have for any \mbox{$j=1,2$}
$$ \tfrac{(\deg Y)^2}{p_j(k)+h} \:\geq\:
\tfrac{\left(2k(d+2-k+r_0)+2h- r_0(r_0+1)\right)^2}{4(p_j(k)+h)} 
=: f_j(k,h,r_0)\,, $$
that is, 
$$\rho_j \:\geq\: \min \, \bigg\{ f_j(k,h,r_0)\;\Big|\;
\begin{array}{l}
1 \leq h \leq \min\,\bigl\{ \tfrac{r_0(r_0+1)}{2},\,\tfrac{k(k-3)}{2}\bigr\}\\
4 \leq k\,, \; 1\leq r_0 \leq k-2 \,.
\end{array}
\bigg\}\,,\qquad j=1,2\,. $$
Remark that for fixed $k,h$ the functions $f_j(k,h,\underline{\phantom{H}}\,)$ 
are increasing in $r_0$ 
(on $\bigl[0,k-\tfrac{1}{2}\bigr]$).  
Hence, they take their minima 
for the minimal possible value, that is, for $r_0$ satisfying 
$$ \left(\,\tfrac{r_0(r_0+1)}{2}=h\,,\;\: r_0\leq k-3\,\right) \;\text{ or }\;
\left(\, r_0=k-2\,,\;\: h \geq \tfrac{(k-3)(k-2)}{2}+1\,\right)\,. $$

\medskip \noindent
{\it Case 1:} \mbox{$r_0=k-2$}, \mbox{$ k(k-3) \geq 2h \geq (k-3)(k-2)+2$}. 
In this case we 
can estimate
$$f_j(k,h,r_0)\:\geq \: \tfrac{\left(2kd + (k-3)(k-2) +2 - 
(k-2)(k-1)\right)^2}{4p_j(k)+2k(k-3)} \: = \: 
\tfrac{2\left(kd -k+3\right)^2}{2p_j(k)+ k(k-3)}\,,$$
whence, due to \mbox{$k\leq (d+3)/2$},
\begin{eqnarray*}
f_1(k,h,r_0) &\geq & \tfrac{(kd-k+3)^2}{k^2-2k+2} \:\geq\:
\tfrac{k^2(d-1)^2}{k^2-2k+2}
\:= \: d^2 + \tfrac{d^2(2k-2)-k^2(2d-1)}{k^2-2k+2} \:\geq \: d^2\,,    \\
f_2(k,h,r_0) &\geq & \left(d-1+\tfrac{3}{k}\right)^2 \:
\geq \: (d-1)^2\,.
\end{eqnarray*}
Thus, (\ref{beta}) is a consequence of  (\ref{cond1}) and (\ref{cond4}). 

\medskip \noindent
{\it Case 2:} \mbox{$h=r_0(r_0+1)/2$}, \mbox{$r_0\leq k-3$}.
It follows that
$$ f_j(k,h,r_0) \:\geq \: \tfrac{2k^2(d+2-k+r_0)^2}{2p_j(k)+r_0(r_0+1)} \,=:\,
g_j(k,r_0)\,,\quad j=1,2\,.
$$ 
We fix \mbox{$k\geq 4$}, and look for the minimum of $g_j(k,r_0)$. Since the 
derivative
$$
\tfrac{\partial}{\partial r_0} \,g_j(k,r_0) \:=\:  
\underbrace{\tfrac{2k^2(d+2-k+r_0)}{(2p_j(k)+r_0(r_0+1))^2}}_{> 0} 
\,\cdot\, 
\bigl(\,r_0(\underbrace{\phantom{\big|}\!\!2k -2d-3}_{<
0})+4p_j(k)+k-d-2\,\bigr) 
$$
changes sign at most once (from positive to negative), the minimum of 
$g_j(k,r_0)$ is taken at one of the endpoints, that is,
$$ \rho_j\:\geq \:\min\, \{ g_j(k,1),\,g_j(k,k-3) \}\,,\quad j=1,2\,.$$ 
We have
\mbox{$\, g_1(k,1)= \tfrac{2k^2}{k^2-k+6} \cdot (d+3-k)^2\,$} and \mbox{$\,
 g_2(k,1)= \tfrac{2k^2}{k^2+3k+2} \cdot (d+3-k)^2$}.
Recall that due to (\ref{Bedingung 1}) we can estimate
$$ d+3-k \:\geq \: \tfrac{d+3}{2}+\sqrt{\tfrac{(d+3)^2}{4} - \deg
 Y}\,,$$ 
whence we obtain
$$ g_1(k,1)\:\geq \: 
\left\{
\renewcommand{\arraystretch}{1.3}
\begin{array}{ll}
\displaystyle
\min\,\bigl\{\,\tfrac{16}{9}\!\;(d-1)^2,\;\,
\tfrac{25}{13} \!\;(d-2)^2\,\bigr\} &\text{ if }k\in \{4,5\}\\
\displaystyle
\tfrac{1}{2} \!\:\bigl(\,
d+3+\sqrt{(d+3)^2 - 4\deg Y}\:\bigr)^2 & \text{ if }
k\geq 6
\end{array}
\right.
$$ 
and
$$ g_2(k,1)\: \geq \:
\tfrac{4}{15} \!\:\bigl(\,
d+3+\sqrt{(d+3)^2 - 4\deg Y}\:\bigr)^2 \,.
$$
On the other hand, we have \mbox{$k<(d+3)/2$}, which implies that
\begin{equation}
\label{beta C}
g_1(k,k-3)\:=\: \tfrac{k^2 (d-1)^2}{k^2-3k+5} \:>\: d^2\,\:\text{ and } \:\,
g_2(k,k-3)\:=\: \tfrac{k^2 (d-1)^2}{k^2-k+3} \:> \: (d-1)^2.
\end{equation}
Thus, if (\ref{cond1}) and (\ref{cond4}) are satisfied and \mbox{$d\geq 8$},
the condition (\ref{beta}) holds whenever 
\begin{equation}
\label{new cond 1}
\textstyle{\sum\limits_{i=1}^s} (\deg Y_i)^2
\:<\: \tfrac{1}{2} \!\:\bigl(\,d+3+\sqrt{(d+3)^2 - 4\deg Y}\:\bigr)^2 
\end{equation}
and
\begin{equation}
\label{new cond 2}
\textstyle{\sum\limits_{z_i \not\in \text{Sing}\,C_k}}
(\deg Y_i)^2 + \textstyle{\sum\limits_{z_i \in \text{Sing}\,C_k}}
\tfrac{(\deg Y_i)^2}{2}\:<\: 
\tfrac{4}{15} \bigl(\,
d+3+\sqrt{(d+3)^2 - 4\deg Y}\:\bigr)^2 .
\end{equation}

\smallskip \noindent
{\it Step 3.} 
In the following, we analyse the conditions (\ref{new cond 1}) and 
(\ref{new cond 2}). We write
\mbox{$\sum \frac{(\deg Y_i)^2}{\varepsilon_i}$} to denote 
the left-hand side of  (\ref{new cond 1}) respectively (\ref{new cond 2}).
As above, we introduce the numbers
\begin{equation}
\label{alpha}
\alpha_{Y,\varepsilon}\::= 
\:\frac{\sum_{i=1}^s \frac{(\deg Y_i)^2}{\varepsilon_i} }{(d+3)^2}
\,, \quad \beta_{Y,\varepsilon}\::=\:\frac{\sum_{i=1}^s 
\frac{(\deg Y_i)^2}{\varepsilon_i} }{\deg  Y} 
\end{equation}
and look for the possible values of $\alpha_{Y,\varepsilon}$ such that 
(\ref{new cond 1}), respectively (\ref{new cond 2}), holds. 
This is the case whenever
$$ \alpha_{Y,\varepsilon} \:< \: \tfrac{K}{4} \cdot \left(\,1+\sqrt{1-
    4\tfrac{\alpha_{Y,\varepsilon}}{\beta_{Y,\varepsilon}}}\:\right)^2,$$
where \mbox{$K=2$}, respectively \mbox{$K=16/15$},
that is, if 
$$\textstyle{\sum\limits_{i=1}^s} \tfrac{(\deg Y_i)^2}{\varepsilon_i}
\:=\: \alpha_{Y,\varepsilon} \cdot (d+3)^2
\:<\:
\tfrac{K\cdot \beta_{Y,\varepsilon}^2}{(\beta_{Y,\varepsilon}+
K)^2} \cdot (d+3)^2\,. $$
Note that this restriction can be reformulated as
$$
\left(\textstyle{\sum_{i=1}^s} \tfrac{(\deg Y_i)^2}{\varepsilon_i}
\right) \cdot
\left( 1+\tfrac{K}{\beta_{Y,\varepsilon}} \right)^2\;< \;K(d+3)^2,
$$
where, by the Cauchy inequality, the left-hand side can be estimated as
\begin{eqnarray*}
\left(\textstyle{\sum\limits_{i=1}^s} \tfrac{(\deg Y_i)^2}{\varepsilon_i}
\right) \cdot
\Bigl( 1+\tfrac{K}{\beta_{Y,\varepsilon}} \Bigr)^2
& = & \frac{ \left(\sum_{i=1}^s \frac{(\deg Y_i)^2}{\varepsilon_i}
+ K\cdot \sum_{i=1}^s \deg Y_i \right)^2}{\sum_{i=1}^s \frac{(\deg
Y_i)^2}{\varepsilon_i}}\\
&\leq & \textstyle{\sum\limits_{i=1}^r}  
\tfrac{\left(\deg Y_i+ K\cdot \varepsilon_i\right)^2}{\varepsilon_i}\,.
\end{eqnarray*}
Finally, since \mbox{$Y_i\subset X'_{\text{fix}}(C,z_i)$}, the conditions
(\ref{new cond 1}) and 
(\ref{new cond 2}) are satisfied if we suppose (\ref{cond2}) and (\ref{cond5}).
\hfill $\Box$

%


\end{document}